

\input epsf.tex

\def\2{{1\over 2}}

\def\d{\delta}
\def\a{\alpha}
\def\b{\beta}
\def\g{\gamma}

\def\s{\sigma}
\def\e{\epsilon}
\def\l{\lambda}

\def\D{\Delta}

\def\fun#1#2#3{#1\colon #2\rightarrow #3}

\def\frac#1#2{{{#1} \over {#2}}}

\def\st{\;\colon\;}
\def\tends{\rightarrow}

\def\dx{\hbox{{\rm d}$x$}}

\def\dr{ {\rm d} }
\def\dt{\hbox{{\rm d}$t$}}

\def\R{{\bf R}}
\def\N{{\bf N}}
\def\Z{{\bf Z}}

\def\T{{\bf T}}

\def\thm#1{\vskip 1 pc\noindent{\bf Theorem #1.\quad}\sl}
\def\lem#1{\vskip 1 pc\noindent{\bf Lemma #1.\quad}\sl}
\def\prop#1{\vskip 1 pc\noindent{\bf Proposition #1.\quad}\sl}
\def\cor#1{\vskip 1 pc\noindent{\bf Corollary #1.\quad}\sl}

\def\proof{\rm\vskip 1 pc\noindent{\bf Proof.\quad}}
\def\fin{\par\hfill $\backslash\backslash\backslash$\vskip 1 pc}
\def\txt#1{\quad\hbox{#1}\quad}
\def\m{\mu}
\def\L{{\cal L}}

\def\s{\sigma}

\def\r{\rho}

\def\2{\frac{1}{2}}
\def\inn#1#2{{\langle #1 ,#2\rangle}}

\def\Mt{{ {\cal M}_1({\bf T}^p)    }}
\def\Co{{{\cal C}_+}}

\def\supp{{\rm supp}}
\def\part{{\partial_{x}}}



\baselineskip= 17.2pt plus 0.6pt
\font\titlefont=cmr17
\centerline{\titlefont Viscous Aubry-Mather theory}
\vskip 1 pc
\centerline{\titlefont and the Vlasov equation}
\vskip 4pc
\font\titlefont=cmr12
\centerline{         \titlefont {Ugo Bessi}\footnote*{{\rm 
Dipartimento di Matematica, Universit\`a\ Roma Tre, Largo S. 
Leonardo Murialdo, 00146 Roma, Italy.}}   }{}\footnote{}{
{{\tt email:} {\tt bessi@matrm3.mat.uniroma3.it}Work partially supported by the PRIN2009 grant "Critical Point Theory and Perturbative Methods for Nonlinear Differential Equations}} 
\vskip 0.5 pc
 
\par
\vskip 2pc
\centerline{\bf Abstract}

The Vlasov equation models a group of particles moving under a potential $V$; moreover, each particle exerts a force, of potential $W$, on the other ones. We shall suppose that these particles move on the $p$-dimensional torus $\T^p$ and that the interaction potential $W$ is smooth. We are going to perturb this equation by a Brownian motion on $\T^p$; adapting to the viscous case methods of Gangbo, Nguyen, Tudorascu and Gomes, we study the existence of periodic solutions and the asymptotics of the Hopf-Lax semigroup.

\vskip 2 pc
\centerline{\bf  Introduction}
\vskip 1 pc

The Vlasov equation models the motion of a group of particles under the action of a time-dependent potential $V$ and a mutual interaction $W$. For definiteness, we shall suppose that the particles move on the torus $\T^p\colon=\frac{\R^p}{\Z^p}$; we put on the position and velocity space $\T^p\times\R^p$ the coordinates $(x,v)$ and we suppose that, at time $t$, the particles are distributed on $\T^p\times\R^p$ according to a probability measure $f_t$. Then, the Vlasov equation has the form
$$\partial_t f_t+\inn{v}{\partial_x f_t}=
{\rm div}_v(f_t\partial_x P_t)   \eqno (VL)_\infty$$
where
$$P(t,x)=V(t,x)+\int_{\T^p\times\R^p}W(x-x^\prime)
\dr f_t(x^\prime,v^\prime)$$
and $\inn{\cdot}{\cdot}$ denotes the scalar product in $\R^p$. Since [7], one looks for weak solutions of $(VL)_\infty$; in other words, given an initial distribution $f_0$, one looks for a continuous curve of probability measures $f_t$ satisfying
$$\int_{\T^p\times\R^p}\phi(0,x,v)\dr f_0(x,v)+$$
$$\int_0^{+\infty}\int_{\T^p\times\R^p}[
\partial_t\phi(t,x,v)-\inn{v}{\partial_x\phi(t,x,v)}+
\inn{\partial_xP(t,x)}{\partial_v\phi(t,x,v)}
]\dr f_t(x,v)  =   0  $$
for all $\phi\in C^\infty_0([0,+\infty)\times\T^p\times\R^p)$.

Our hypotheses on $V$ and $W$ are the following:

\noindent 1) $V\in C(\T,C^3(\T^p))$, and

\noindent 2) $W\in C^3(\T^p)$. Thus, $W$ lifts to a $C^3$ function on $\R^p$, $\Z^p$-periodic; we shall also suppose that 
$W(x)=W(-x)$ and that $W(0)=0$.

\vskip 1pc

A recent idea (see [11], [12]) is to view $(VL)_\infty$ as a Lagrangian system in the space of measures; indeed, it is possible to define what it means for a curve $\mu_t$ of probability measures on $\T^p$ to minimize the Lagrangian action
$$\int_{t_0}^{t_1}\left[
\2||\dot\mu_t||^2+
\int_{\T^p}V(t,x)\dr\mu_t+
\2\int_{\T^p\times\T^p}W(x-x^\prime)\dr(\mu_t\times\mu_t)
(x,x^\prime)
\right]\dt  .  \eqno (1)$$
The advantages are that one can use the tools of Lagrangian dynamics (Aubry-Mather theory, Hamilton-Jacobi equations, minimal characteristics, etc...) albeit in the difficult "differential manifold" of probability measures. 

In this paper, we are going to adapt to the viscous case an older approach: following [7], we jury-rig a fixed point theorem to the viscous Mather theory of [13]. Let us briefly outline what we are doing in the case of  periodic orbits.

We let $\psi(t,x)$ be a continuous family of densities, periodic in time; in other words, we ask that

\noindent d1) $\psi\in C(\T\times\T^p)$

\noindent d2) $\psi\ge 0$

\noindent d3) $\int_{\T^p}\psi(t,x)\dx=1$ for all $t\in\T$.

Let us define
$$P_\psi(t,x)=V(t,x)+\int_{\T^p}W(x-x^\prime)\psi(t,x^\prime)
\dr x^\prime  $$
and, for $c\in\R^p$, let us set
$$\fun{\L_{c,\psi}}{\T\times\T^p\times\R^p}{\R},\qquad
\L_{c,\psi}(t,x,\dot x)=\2|\dot x|^2-\inn{c}{\dot x}-P_\psi(t,x)$$
$$\fun{H_{\psi}}{\T\times\T^p\times\R^p}{\R},\qquad
H_\psi(t,x,p)=\2|p|^2+P_\psi(t,x)  .  $$
We have the following.

\thm{1} Let $c\in\R^p$ and let $\b>0$. Then, there is a couple of functions
$$\rho_\b\in C^1(\T\times\T^p)\cap C(\T,C^2(\T^p)),
\qquad
u_\b\in C^1(\T,C^1(\T^p))\cap C(\T,C^3(\T^p))$$
and $\bar H_{\b}(c)\in\R$ such that $\rho_\b$ satisfies points d1)-d3) above and
$$\frac{1}{2\b}\D u_\b+\partial_t u_\b -
H_{\rho_\b}(t,x,c-\partial_xu_\b)+\bar H_{\b}(c)=0  ,
\eqno (HJ)_{\rho_\b,per}$$
$$\frac{1}{2\b}\D\r_\b-{\rm div}[\rho_\b\cdot(c-\partial_xu_\b)]-
\partial_t\rho_\b=0  .   \eqno(FP)_{c-\partial_x u_\b,per}  $$
Moreover, among the triples $(\rho_\b,u_\b,\bar H_{\b}(c))$ which solve $(HJ)_{\rho_\b,per}-(FP)_{c-\partial_x u_\b, per}$, there is one which minimizes 
$$\int_{\T\times\T^p}
\L_{c,\2\rho_\b}(t,x,\partial_x u_\b)\rho_{\b}(t,x)\dr t\dx  .  
\eqno (2)$$

\rm

\vskip 1pc

Thus, our "characteristics" are the solutions of a 
Fokker-Planck equation bringing mass forward in time; the drift of this equation, or the optimal trajectory, is determined by a Hamilton-Jacobi equation, backward in time. This is quite typical for this kind of problems: see for instance equation (5.40) of [17] or theorem 3.9 of [12]. 

We briefly sketch the proof of theorem 1; the complete details are in section 1 below. First of all, we fix $\psi$ satisfying points 
d1)-d3) above; then we find, as in [13], a couple 
$(u_\psi,\bar H_\psi(c))$ which solves $(HJ)_{\psi,per}$. By [13], the number $\bar H_\psi(c)$ is unique and $u_\psi$ is unique  up to an additive constant. To $c-\partial_x u_\psi$ is associated a stochastic flow, whose stationary Fokker-Planck equation is 
$(FP)_{c-\partial_xu_\psi,per}$; again by [13], 
$(FP)_{c-\partial_xu_\psi,per}$ has a unique periodic solution 
$\rho_\psi$ satisfying d1)-d3). In other words, we have a map 
$\fun{}{\psi}{\rho_\psi}$ bringing densities to densities; we shall find a fixed point $\rho_\b$ of this map by the Schauder fixed point theorem. We shall see that 
$(u_{\rho_\b},\rho_\b,\bar H_{\rho_\b}(c))$ solves 
$(HJ)_{\rho_\b,per}-(FP)_{c-\partial_x u_{\rho_\b}, per}$ practically by definition; the existence of a minimum in (2) will follow from the fact that the fixed points of $\rho_\b$ are a compact set.

In section 2, we study the Hopf-Lax semigroup. We denote by 
$\Mt$ the space of Borel probability measures on $\T^p$ with the 
1-Wasserstein distance (see section 1 for a definition); we shall prove the following.

\thm{2} Let $\fun{U}{\Mt}{\R}$ be of the form
$$U(\mu)=\int_{\T^p}f\dr\mu$$
for some $f\in C^3(\T^p)$. Let $\mu\in\Mt$ and let $m\in\N$. Then, there are $R_\b\in C([-m,0],\Mt)$ with density 
$\r_\b\in C^1((-m,0]\times\T^p)\cap C((-m,0]\times C^2(\T^p))$ and 
$u_\b\in C^1([-m,0],C^1(\T^p))\cap C([-m,0],C^3(\T^p))$ such that $u_\b$ solves
$$\left\{
\eqalign{
{}&\frac{1}{2\b}\D u_\b+\partial_tu_\b-
H_{\rho_\b}(t,x,c-\partial_xu_\b)=0,\qquad t<0\cr
{}& u_\b(0,x)=f\quad\forall x\in\T^p     
}      \right.       \eqno (HJ)_{\rho_\b,f}$$
and $R_\b$ together with its density $\r_\b$ solve
$$\left\{
\eqalign{
{}&\frac{1}{2\b}\D \r_\b-{\rm div}[\r_\b\cdot(c-\partial_x u_\b)]-
\partial_t\rho_\b=0, \qquad t>-m  \cr
{}&R_\b(-m)=\mu   .   
}       \right.     \eqno (FP)_{-m,c-\partial_xu_\b,\mu}    $$
Among the solutions $(u_\b,\r_\b)$ of 
$(HJ)_{\b,f}$-$(FP)_{-m,c-\partial_xu_\b,\mu}$, there is one which minimizes
$$\int_{-m}^0\dt\int_{\T^p}\L_{c,\2\rho_\b}(t,x,c-\partial_xu_\b)
\rho_\b\dx+U(\rho_\b(0))  .  $$
We call such a minimum $(\Lambda^m_c U)(\mu)$.

\rm
\vskip 1pc

Since minimizing over fixed points is uncomfortable, one could ask whether this restriction can be removed, getting a problem more similar to (1).

\thm{3} 1) Let $U$ and $(\Lambda^m_c U)(\mu)$ be as in theorem 2. Then,
$$(\Lambda^m_cU)(\mu)=\min_Y E_w\left\{
\int_{-m}^0\L_{c,\2\r}(t,X,Y)\dt\right\}         
+U(\r(0))     \eqno (3)$$
where $\r$ solves $(FP)_{-m,Y,\mu}$, $E_w$ denotes expectation with respect to the Brownian motion $w$, $X$ solves the stochastic differential equation
$$\left\{
\eqalign{
{}&\dr X(-m,s,x)=Y(s,X(-m,s,x))\dr t+\dr w(s)\quad 
s\ge -m\cr
{}&X(-m,-m,x)=X_\mu    
}
\right.     \eqno (SDE)_{-m,Y,\mu}$$
for a random variable $X_\mu$ of law $\mu$, independent on 
$w(s)$ for $s\ge-m$. The minimum is taken over all Lipschitz vector fields $Y$.

\noindent 2) Any minimal $Y$ satisfies $Y=c-\partial_x u$, where $u$ solves $(HJ)_{\r,f}$.

\vskip 1pc

\rm

In other words, $(HJ)_{\r_\b,f}$-$(FP)_{-m,c-\partial_xu_\b,\mu}$ are the Euler-Lagrange equations of the functional (3), exactly as in the zero-viscosity situation (we refer again to [12], theorem 3.9.) We also note a quirk of the notation: in the Hamilton-Jacobi equation we have $H_{\r_\b}$, while in (3) we have 
$\L_{c,\2\r_\b}$; again, we share this factor two with the 
zero-viscosity situation and we shall see the reason for it in the proof of lemma 2.5 below.

Theorems like theorem 3 are common in the theory of mean field games (see for instance [4]). In the language of mean field games, we are saying that each particle tries to minimize unilaterally the cost 
$$\min_Y E_w\left\{
\int_{-m}^0\L_{c,\r}(t,X,Y)\dt+f(X(0))
\right\}  $$
where $X$ solves $(SDE)_{-m,Y,\d_{x_0}}$ and $\r$ is the distribution of the other particles; this is the reason of equation 
$(HJ)_{\r_\b,f}$ in theorem 2. The result of the independent efforts of all the particles (or the Nash equilibrium, as it is called) is that the whole community minimizes (3).

Let $\fun{U}{\Mt}{\R}$ be bounded; theorem 3 prompts us to define
$$(\Psi^m_cU)(\mu)=\inf_Y \left\{
E_w\int_{-m}^0\L_{c,\2\r}(t,X,Y)\dt+U(\r(0))
\right\}  \eqno (4)$$
where the infimum is taken over all Lipschitz vector fields $Y$; the density $\r$ satisfies $(FP)_{-m,Y,\mu}$. Naturally, if $U$ is linear as in theorem 3, then $\Psi^m_cU=\Lambda^m_cU$.

We shall see in proposition 2.10 below that $\Psi^m_c$ has the semigroup property: 
$\Psi^{m+n}_c=\Psi^m_c\circ\Psi^n_c$.

Theorem 3 tells us that the infimum in (4) is a minimum when $U$ is a linear function on measures as in theorem 2; we don't know whether this is true when $U$ is in some more reasonable class, for instance continuous or Lipschitz. We don't even know whether, for $U$ continuous, 
$\Psi_c^1U$ is continuous; however, when $U$ is linear as in theorem 2, we can prove that $\Psi^m_cU$ is Lipschitz, uniformly in $m$. This allows us to find, for a suitable $\l\in\R$, Lipschitz fixed points of the operator $\Psi_{c,\l}$ defined by
$$\fun{\Psi_{c,\l}}{U}{\Psi_c^1U+\l}  .  $$

\thm{4} There is a unique $\l\in\R$ for which $\Psi_{c,\l}$ has a fixed point $\hat U$ in $C(\Mt,\R)$. In other words, for any 
$\mu\in{\cal M}_1(\T^p)$, there is a Lipschitz vector field $\bar Y$ such that
$$(\Psi^1_{c,\l}\hat U)(\mu)=
E_w\left\{
\int_{-1}^0
\L_{c,\2\bar\r}(t,\bar X,\bar Y)\dt
\right\}     +\hat U(\bar\r(0)),    \eqno(5)$$
where $X$ solves $(SDE)_{-1,\bar Y,\mu}$ and $\bar\r$ solves
$(FP)_{-1,\bar Y,\mu}$.

The function $\hat U$ is Lipschitz for the 1-Wasserstein distance; by (5), the infimum in the definition of (4) of $\Psi^1_{c,\l}\hat U$ is a minimum.

\rm
\vskip 1pc

The proof of this theorem is similar to the corresponding statement in Aubry-Mather theory. Indeed, using an approximation with finitely many particles, we shall prove that, for a suitable $\l\in\R$, the sequence $(\Lambda_{c,\l})^n(0)$ of continuous functions on the compact space $\Mt$ is equibounded and equilipschitz; by
Ascoli-Arzel\`a, it has a subsequence converging to a limit 
$\hat U$; we shall prove that $\hat U$ is a fixed point of 
$\Lambda_{c,\l}$.

\vskip 2pc
\centerline{\bf \S 1}
\centerline{\bf Periodic orbits}
\vskip 1pc

In this section, we are going to prove theorem 1. We begin with a study of $(HJ)_{\psi,per}$; we follow the approach of [13] but, for completeness' sake, we reprove several results of this paper using, as in [2], the Feynman-Kac formula. 

\vskip 1pc

\noindent{\bf Definitions.} 

\noindent $\bullet$ We group in a set $Den$ the functions on $\T\times\T^p$ which satisfy points d1)-d3) in the introduction. Clearly, the set $Den$ is closed in $C(\T\times\T^p)$.

\vskip 1pc

\noindent$\bullet$ We define $\Mt$ as the space of all Borel probability measures on $\T^p$; if $\mu_1,\mu_2\in\Mt$, we define the $1$-Wasserstein distance between them as 
$$d_1(\mu_1,\mu_2)=\min\{
\int_{\T^p\times\T^p}|x-x^\prime|_{\T^p}\dr\g(x,x^\prime)
\}$$
where $|x-x^\prime|_{\T^p}$ is the distance on the flat torus 
$\T^p$. The minimum is taken over all the Borel probability measures $\g$ on $\T^p\times\T^p$ whose first and second marginals are, respectively, $\mu_1$ and $\mu_2$. It is standard (see for instance section 7.1 of [17]) that $d_1$ turns $\Mt$ into a complete metric space, and induces the weak$\ast$ topology.

We note that, if $\psi\in Den$ and $\L^p$ denotes the Lebesgue measure on $\T^p$, then the function $\fun{}{t}{\psi(t,\cdot)\L^p}$ belongs to $C(\T,\Mt)$.

\vskip 1pc

\noindent$\bullet$ We extend the definition of $P_\psi$ we gave in the introduction: for $\psi\in C(\R,\Mt)$ we set
$$P_\psi(t,x)=V(t,x)+\int_{\T^p}
W(x-x^\prime)\dr\psi(t,x^\prime)
   \eqno (1.1) $$

\lem{1.1} There is $C_1>0$, independent on $\psi\in C(\R,\Mt)$, such that the function $P_\psi(t,x)$ defined in (1.1) satisfies 
$$||P_\psi||_{C(\R,C^3(\T^p))}\le C_1  .  \eqno (1.2)$$

\proof We recall that, by definition, 
$$|| P_\psi ||_{C(\R,C^3(\T^p))}=
\sup_{t\in\R}||
P_\psi(t,\cdot)
||_{C^3(\T^p)}     $$
where, as usual, 
$$||f||_{C^3(\T^p)}=||f||_{C^0(\T^p)}+||D_x f||_{C^0(\T^p)}+
||D^2_x f||_{C^0(\T^p)} +||D^3_x f||_{C^0(\T^p)}   .  $$
By our hypotheses on $V$ and $W$, we have that
$$|| V ||_{C(\R,C^3(\T^p))}+|| W ||_{C^3(\T^p)}=C_1<+\infty  .
\eqno (1.3)  $$
For $0\le j\le 3$, differentiation under the integral sign implies that
$$D^j_x P_\psi(t,x)=D^j_x V(t,x)+
\int_{\T^p}D^j_x W(x-x^\prime)\dr\psi(t,x^\prime)  .  $$
Since $\fun{}{t}{\psi(t,\cdot)}$ is continuous from $\R$ to the weak$\ast$ topology, the formula above implies that $P_\psi$ is in 
$C(\R,C^3(\T^p))$. Since $\psi(t,\cdot)$ is a probability mesure and the $C^3$ norm is convex, (1.2) follows from the last formula and (1.3).

\fin

From now on, we shall fix $\psi\in Den$; the functions $P_\psi$, 
$\L_{c,\psi}$ and $H_\psi$ are defined as in the introduction. Following [2], we note that, if $(u,A)$ solves
$$\left\{
\eqalign{
{}&\frac{1}{2\b}\Delta u+\partial_t u-H_\psi(t,x,c-\partial_x u)+A=0
\quad\forall t\in\R\cr
{}&u(t,\cdot)=u(t+1,\cdot)\quad\forall t\in\R
}   \right.  
\eqno (HJ)_{\psi,per}$$ 
and is periodic in space (i. e. it quotients to a continuous function on $\T^p$), then the couple 
$(v,A)=(e^{-\b u},A)$ is a solution, periodic in space, of the  "twisted" Schr\"odinger equation
$$\left\{
\eqalign{
{}&\partial_tv+e^{-\b\inn{c}{x}}\left[
\frac{1}{2\b}\D+\b P_\psi(t,x)-\b A
\right]
(e^{\b\inn{c}{x}}v)=0
\quad\forall t\in\R\cr
{}&v(t,\cdot)=v(t+1,\cdot)\quad\forall t\in\R   .    }
\right.      \eqno (TS)_{\psi,per}   $$
Vice-versa, the logarithm of a positive solution of $(TS)_{\psi,per}$ solves $(HJ)_{\psi,per}$. Thus, solving $(HJ)_{\psi,per}$ reduces to solving $(TS)_{\psi,per}$; that's what we are going to do next.

For $\phi\in C(\T^p)$ and $A\in\R$ we consider the evolution (or involution, since it goes backward in time) problem
$$\left\{
\eqalign{
{}&\partial_tv+e^{-\b\inn{c}{x}}\left[
\frac{1}{2\b}\D+\b P_\psi(t,x)-\b A
\right]
(e^{\b\inn{c}{x}}v)=0, \quad t\le 0\cr
{}&v(0,x)=\phi(x)   .    }
\right.      \eqno (TS)_{\psi,\phi}   $$
If $t\le 0$, we can use the Feynman-Kac formula (see for instance [6]) and write the unique solution of $(TS)_{\psi,\phi}$ as
$$v(t,x)=(L_{(\psi,A,t)}\phi)(x),\qquad t\le 0$$
where
$$(L_{(\psi,A,t)}\phi)(x)=e^{-\b\inn{c}{x}}\cdot
E_w\left[
e^{
\int_t^0\b[
P_\psi(\tau,\frac{1}{\sqrt\b}w(\tau)+x)-A
]\dr\tau}
e^{
\b\inn{c}{\frac{1}{\sqrt\b}w(0)+x}
}
\phi\left(
\frac{1}{\sqrt\b} w(0)+x
\right)
\right]    .    \eqno (1.5)$$
In the formula above, $w$ is a Brownian motion on $[t,+\infty]$ with $w(t)=0$, and $E_w$ is the expectation with respect to the Wiener measure.

We shall see in lemma 1.4 below that there is a bijection between the positive eigenfunctions of $L_{(\psi,0,-1)}$ and the positive solutions of $(TS)_{\b,per}$; now, we prove that such eigenfunctions exist.

\lem{1.2} 1) (Existence) there is $(v,B)\in C(\T^p)\times\R$ such that
$$\left\{
\eqalign{
{}&L_{(\psi,0,-1)}v=Bv\cr
{}&v>0\cr
{}&B>0  .   
}
\right.    \eqno (1.6)_\psi  $$

\noindent 2) (Uniqueness) Let $(v_1,B_1)$ and $(v_2,B_2)$ solve 
$(1.6)_\psi$; then, $B_1=B_2$ and $v_1=\a v_2$ for some 
$\a>0$. In particular, there is a unique couple $(v_\psi,B_\psi)$ which satisfies $(1.6)_\psi$ and such that
$$\int_{\T^p}v_\psi(x)\dx=1  .  \eqno (1.7)  $$

\proof We recall from [16] (see also chapter XVI of [3] for G. Birkhoff's original exposition) a few facts about the Perron-Frobenius theorem. Let us denote by $\Co\subset C(\T^p)$ the cone of strictly positive, continuous functions. We forego the easy proof that $L_{(\psi,0,-1)}$ brings $\Co$ into itself.

Let $v_1,v_2\in\Co$; we say that $v_1$ and $v_2$ are equivalent, or $v_1\simeq v_2$, if $v_1=t v_2$ for some $t>0$. Given $v_1,v_2\in\Co$, we define
$$\a(v_1,v_2)=\sup\{
t>0\st v_2-tv_1\in\Co
\}  $$
and 
$$\theta(v_1,v_2)=-\log[
\a(v_1,v_2)\a(v_2,v_1)
]  .  $$
It turns out ([16]) that  $(\frac{\Co}{\simeq},\theta)$ is a complete metric space. We refer again to [16] or [3] for the proof that
$$\theta(
L_{(\psi,0,-1)}v_1,L_{(\psi,0,-1)}v_2
)  \le   (1-e^{-D})\theta(v_1,v_2)$$
where
$$D=\sup_{v_1,v_2\in\Co}
\theta(
L_{(\psi,0,-1)}v_1,L_{(\psi,0,-1)}v_2
)  .  $$
As a consequence, points 1) and 2) follow from the contraction mapping theorem if we prove that $D<+\infty$. Actually, we are going to show that $D$ is bounded from above independently of 
$\psi\in Den$; equivalently, the Lipschitz constant of 
$L_{(\psi,0,-1)}$ does not depend on $\psi$. We shall need this fact in the next lemma.

Let $v_1,v_2\in\Co$. Recalling the definition of $\theta$, we see that
$$\theta(v_1,v_2)\le\log\left(
\frac{\max v_2}{\min v_1}\cdot\frac{\max v_1}{\min v_2}
\right)  .  \eqno (1.8)$$
Thus, $D<+\infty$ follows if we prove that there is $C_3>0$ such that
$$\frac{
\max L_{(\psi,0,-1)}v
}{
\min L_{(\psi,0,-1)}v
}    \le C_3     \eqno (1.9)  $$
for all $v\in\Co$; since the term on the left is homogeneous of degree zero in $v$, we can suppose that $v$ satisfies (1.7). 

We prove (1.9); in the following, $C_i$ always denotes a constant independent on $v$ and $\psi$. By (1.5) and the fact that $v>0$, we have that
$$(
L_{(\psi,0,-1)}v
)(x)\ge
e^{
-\b\inn{c}{x}
}   
e^{
\b\min P_\psi
}
\frac{1}{\sqrt{(2\pi)^p}}
\int_{\R^p}
e^{
\b\inn{c}{\frac{1}{\sqrt\b}z+x}
}
v\left(
\frac{1}{\sqrt\b}z+x
\right)
e^{
-\frac{|z|^2}{2}
}    \dr z  .  $$
Setting $\frac{1}{\sqrt\b}z=y$ and simplifying $e^{-\b\inn{c}{x}}$ outside the integral with $e^{\b\inn{c}{x}}$ inside, we get the first inequality below
$$(
L_{(\psi,0,-1)}v
)(x)\ge
e^{
\b\min P_\psi
}
\sqrt{
\left(
\frac{\b}{2\pi}
\right)^p
}
\int_{\R^p}
e^{
\b\inn{c}{y}
}
v(x+y)
e^{
-\frac{\b}{2}|y|^2
}
\dr y\ge$$
$$e^{-\b C_1}
\sqrt{
\left(
\frac{\b}{2\pi}
\right)^p
}
\int_{[0,1]^p}e^{
\b\inn{c}{y}
}
v(x+y)
e^{
-\frac{\b}{2}|y|^2
}
\dr y  .   $$
The second inequality above comes from lemma 1.1 and the fact that $v$, which belongs to $\Co$, is positive. By lemma 1.1, the constant $C_1$ does not depend on $\psi\in Den$.

We assert that
$$(
L_{(\psi,0,-1)}v
)(x)\ge
e^{-\b C_1}
\sqrt{
\left(
\frac{\b}{2\pi}
\right)^p
}
\min_{y\in[0,1]^p}
\left[
e^{\b\inn{c}{y}}e^{-\frac{\b}{2}|y|^2}
\right]  \int_{
[0,1]^p
}  
v(x+y)\dr y= C_5    \eqno (1.10)$$
for a constant $C_5>0$ independent on $\psi$ and $v$. Indeed, the inequality follows since $v$ is positive; since $v$ is periodic and satisfies (1.7), the integral above is $1$, and the equality follows.

For the estimate from above, we get from (1.5) that
$$(
L_{(\psi,0,-1)}v
)(x)\le
e^{
-\b\inn{c}{x}
}
e^{
\b\max P_\psi
}
\frac{1}{\sqrt{(2\pi)^p}}\int_{\R^p}
e^{
\b\inn{c}{\frac{1}{\sqrt\b}z+x}
}
v(\frac{1}{\sqrt\b}z+x)
e^{
-\frac{|z|^2}{2}
}
\dr z  .  $$
We simplify $e^{-\b\inn{c}{x}}$ outside the integral with 
$e^{\b\inn{c}{x}}$ inside; now lemma 1.1 gives us the first inequality below; the equality follows from the change of variables 
$\frac{1}{\sqrt\b}z=y$.
$$(
L_{(\psi,0,-1)}v
)(x)\le
\frac{e^{\b C_1}}{\sqrt{(2\pi)^p}}\int_{\R^p}
e^{
\b\inn{c}{\frac{1}{\sqrt\b}z}
}
v(\frac{1}{\sqrt\b}z+x)
e^{
-\frac{|z|^2}{2}
}   \dr z  =$$
$$e^{\b C_1}
\sqrt{
\left(
\frac{\b}{2\pi}
\right)^p
}
\int_{\R^p}
e^{\b\inn{c}{y}}
v(x+y)
e^{
-\frac{\b}{2}|y|^2
}
\dr y  .  $$
Since $v$ is positive periodic, and by (1.7) integrates to $1$ on the unit cube, we get the first inequality below.
$$(
L_{(\psi,0,-1)}v
)(x)\le
e^{\b C_1}
\sqrt{
\left(
\frac{\b}{2\pi}
\right)^p
}\cdot
\sum_{k\in\Z^p}\max_{y\in k+[0,1]^p}
\left[
e^{
\b\inn{c}{y}
}
e^{
-\frac{\b}{2}|y|^2
}
\right]  \le$$
$$C_6\sum_{k\in\Z^p}
e^{
\b|c|(|k|+\sqrt p)
}
e^{
-\frac{\b}{2}(|k|-\sqrt{p})^2
}   .  $$
Since the sum in the last formula is finite, we get that
$$(
L_{(\psi,0,-1)}v
)(x)\le C_7\qquad\forall x\in\T^p     \eqno (1.11)$$
for a constant $C_7>0$ independent on $\psi$ and $v$.
Now (1.9) follows from the last formula and (1.10); we have seen that (1.9), by the contraction mapping theorem, implies points 1) and 2) of the thesis.

\fin

\lem{1.3} Let $\psi\in Den$ and let $(v_\psi,B_\psi)$ be as in the last lemma. Then, 
$v_\psi\in C^3(\T^p)$ and the following two points hold.

\noindent 1) (Uniform estimates) There is $C_8>0$, independent on $\psi\in Den$, such that
$$||v_\psi||_{C^3(\T^p)}\le C_8,\leqno i) $$
$$\frac{1}{C_8}\le v_\psi(x)\le C_8\quad\forall x\in\T^p
\leqno ii)$$ 
$$\frac{1}{C_8}\le B_\psi\le C_8  .  \leqno iii)  $$

\noindent 2) (Continuous dependence) The function 
$$\fun{K}{Den}{C^3(\T^p)\times\R},\qquad
\fun{K}{\psi}{(v_\psi,B_\psi)}$$
is continuous.

\proof We prove point 1). Since $v_\psi$ satisfies (1.7), by (1.10) and (1.11) there is $C_8>1$ such that
$$\frac{1}{C_8}\le\min L_{(\psi,0,-1)}v_\psi\le
\max L_{(\psi,0,-1)}v_\psi\le C_8  .  $$
Since we also have that $L_{(\psi,0,-1)}v_\psi=B_\psi v_\psi$, we get that
$$\frac{1}{B_\psi}\cdot\frac{1}{C_8}\le\min v_\psi\le v_\psi\le
\max v_\psi\le\frac{1}{B_\psi} C_8  .  \eqno (1.12)$$
Integrating on $\T^p$ and using (1.7), we get that 
$$\frac{1}{B_\psi}\cdot\frac{1}{C_8}\le
\int_{\T^p}v_\psi\dx=1\le\frac{1}{B_\psi}\cdot C_8$$
from which $iii$) of point 1) follows.

From point $iii$) and (1.12), possibly increasing $C_8$, we get point $ii$). We show $i$). 

We would like to differentiate under the integral sign in (1.5); we cannot do this immediately, because we only know that the final condition $\phi$ (which in our case is $v_\psi$) is in $C^0$.
Let $E_{(0,z)}$ denote the expectation of the Brownian bridge with $w(-1)=0$ and $w(0)=z$; by (1.5) we get that, for $v\in\Co$, 
$$(
L_{(\psi,0,-1)}v
)(x)=$$
$$e^{
-\b\inn{c}{x}
}
\frac{1}{\sqrt{(2\pi)^p}}\cdot
\int_{\R^p}
e^{
-\frac{|z|^2}{2}
}
e^{
\b\inn{c}{\frac{1}{\sqrt\b}z+x}
}
v(\frac{1}{\sqrt\b}z+x)\cdot
E_{(0,z)}\left[
e^{
\int_{-1}^0\b P_\psi(\tau,\frac{1}{\sqrt\b}w(\tau)+x)\dr\tau
}
\right]   \dr z  .  $$
Setting $\frac{1}{\sqrt\b}z+x=y$, we get that
$$(
L_{(\psi,0,-1)}v
)(x)=
\sqrt{
\left(
\frac{\b}{2\pi}
\right)^p
}  \cdot e^{-\b\inn{c}{x}}\cdot$$
$$\int_{\R^p}
e^{
-\frac{\b}{2}|y-x|^2
}
e^{
\b\inn{c}{y}
}
v(y)\cdot
E_{(0,\sqrt\b(y-x))}\left[
e^{
\int_{-1}^0\b P_\psi(\tau,\frac{1}{\sqrt\b}w(\tau)+x)\dr\tau
}
\right]   \dr y  .  $$
We recall from [14] that, if $\tilde w$ is a Brownian bridge with 
$\tilde w(-1)=\tilde w(0)=0$, then 
$w(t)\colon=\sqrt\b (y-x)(t+1)+\tilde w(t)$ is a Brownian bridge with $w(-1)=0$, $w(0)=\sqrt\b (y-x)$. This and the last formula imply that
$$(
L_{(\psi,0,-1)}v
)(x)=
\sqrt{
\left(
\frac{\b}{2\pi}
\right)^p
}  \cdot e^{-\b\inn{c}{x}}\cdot$$
$$\int_{\R^p}
e^{
-\frac{\b}{2}|y-x|^2
}
e^{
\b\inn{c}{y}
}
v(y)\cdot
E_{(0,0)}\left[
e^{
\int_{-1}^0\b P_\psi(\tau,(y-x)(\tau+1)+\frac{1}{\sqrt\b}
\tilde w(\tau)+x)\dr\tau
}
\right]   \dr y  .  $$
The formula above allows us to differentiate under the integral sign, even if $v$ is only $C^0$; using lemma 1.1, we easily get
$$||
L_{(\psi,0,-1)}v
||_{C^3(\T^p)}\le C_9||v||_{C^0(\T^p)}     \eqno (1.13)$$
for a constant $C_{9}$ independent of $\psi$. By $ii$), we get that 
$$||
L_{(\psi,0,-1)}v
||_{C^3(\T^p)}\le C_8\cdot C_9  .  $$
Since $L_{(\psi,0,-1)}v_\psi=B_\psi v_\psi$, formula $i$) now follows from $iii$).

We prove point 2); in the first three steps below, we show a weaker result, namely that the map $\fun{}{\psi}{v_\psi}$ is continuous from $Den$ to $C^0(\T^p)$; this will follow from the theorem of contractions depending on a parameter applied to the map
$$\fun{\Xi}{(Den,||\cdot||_{\sup})\times({\cal C}_+,\theta)}{
({\cal C}_+,\theta)}  ,   \qquad
\fun{\Xi}{(\psi,v)}{L_{(\psi,0,-1)}v}   .    $$

\noindent{\bf Step 1.} We begin to observe that $\theta$ and the 
$\sup$ norm induce equivalent topologies on the subset ${\cal A}$ of the functions of $\Co$ which satisfy (1.7). Indeed, (1.8) proves that the $C^0$ topology is stronger; for the opposite inclusion, let 
$\theta(v_n, v)\tends 0$ and let $v_n$, $v$ satisfy (1.7). Since 
$\theta(v_n, v)\tends 0$, we have that, for any $\e>0$ and $n$ large enough,
$$\frac{1-\e}{\a(v_n,v)}\le\a(v,v_n)\le\frac{1+\e}{\a(v_n,v)}  .  $$
The definition of $\a$ implies the first two inequalities below; the last one follows by the first inequality above.
$$\a(v,v_n)v\le v_n\le\frac{1}{\a(v_n,v)}v\le
\frac{\a(v,v_n)}{1-\e}v  .  $$
Since $v$ and $v_n$ satisfy (1.7), if we integrate the formula above on $\T^p$, we get that 
$\a(v_n,v)\tends 1$ and that $\a(v,v_n)\tends 1$; since 
$\min v>0$, again from the formula above we get that 
$v_n\tends v$ uniformly.

\noindent{\bf Step 2.} Let $v\in\Co$ be fixed; we assert that the map $\fun{}{\psi}{\Xi(\psi,v)}$ is continuous from the 
$||\cdot||_{\sup}$ to the $\theta$ topology. Indeed, we saw in step 1 that, on $\Co$, the $C^0$ topology is stronger than the 
$\theta$ topology; thus, it suffices to prove that 
$\fun{\Xi(\cdot,v)}{(Den,||\cdot||_{\sup})}{(\Co,||\cdot||_{\sup})}$ is continuous. The proof of this, which ends the proof of the assertion, follows by applying the theorem of continuity under the integral sign to (1.5), and we forego it.  

\noindent{\bf Step 3.} We assert that the map 
$\fun{}{\psi}{v_\psi}$ is continuous from $(Den, ||\cdot||_{\sup})$ to 
$({\cal A},||\cdot||_{\sup})$; by step 1, it suffices to prove that it is continuous from $(Den, ||\cdot||_{\sup})$ to 
$({\cal A},\theta)$. We have seen in the proof of lemma 1.2 that 
$\fun{}{v}{\Xi(\psi,v)}$ is a contraction for the $\theta$-topology, whose Lipschitz constant does not depend on $\psi$. Since $\fun{}{\psi}{\Xi(\psi,v)}$ is continuous by step 2, we can apply the theorem of contractions depending on a parameter, and get that the map 
$\fun{}{\psi}{v_\psi}$ is continuous from $(Den, ||\cdot||_{\sup})$ to $({\cal C}_+,\theta)$, as we wanted.

\noindent{\bf Step 4.} We assert that the map  
$\fun{}{\psi}{B_\psi}$ is continuous from $Den$ to $\R$. Since 
$L_{(\psi,0,-1)}v_\psi=B_\psi v_\psi$, it suffices to prove that both maps $\fun{}{\psi}{v_\psi}$ and $\fun{}{\psi}{L_{(\psi,0,-1)}v_\psi}$ are continuous from $Den$ to $C^0(\T^p)$. The first fact has been proven in step 3; we prove that $\fun{}{\psi}{L_{(\psi,0,-1)}v_\psi}$ is continuous. Indeed, 
$$||L_{(\psi^\prime,0,-1)}v_{\psi^\prime}-
L_{(\psi,0,-1)}v_{\psi}||_{\sup}\le
||L_{(\psi^\prime,0,-1)}(v_{\psi^\prime}-v_\psi)||_{\sup}+
||(L_{(\psi^\prime,0,-1)}-L_{(\psi,0,-1)})v_\psi||_{\sup} . $$
Now the assertion follows from the fact that (with the $\sup$ norm in all spaces) $\fun{}{\psi}{v_\psi}$ is continuous, that $\fun{}{\psi}{L_{(\psi,0,-1)}v}$ is continuous, and that $\fun{}{v}{L_{(\psi,0,-1)}v}$ is uniformly Lipschitz by (1.13).

\noindent{\bf End of the proof of point 2).} For $\phi\in\Co$, we get from (1.5) that
$$L_{(\psi,r,-1)}\phi=
e^{-\b r}L_{(\psi,0,-1)}\phi .  \eqno (1.14)$$
Setting $A_\psi=\frac{1}{\b}\log B_\psi$, the formula above implies that
$$L_{(\psi,A_\psi,-1)}v_\psi=v_\psi  .  \eqno (1.15)$$
The same proof which yielded (1.13) also yields that there is 
$C_{10}>0$ such that, if $A$ and $A^\prime$ satisfy the estimate of point 1), $iii$) of this lemma, then
$$||
L_{(\psi,A,-1)}v-L_{(\psi^\prime,A^\prime,-1)}v
||_{C^3(\T^p)}   \le
C_{10}(
||
\psi-\psi^\prime
||_{C^0(\T^p)}+
|A-A^\prime|
)
\cdot ||v||_{C^0(\T^p)}  .  \eqno (1.16)$$
Thus,
$$||
v_\psi-v_{\psi^\prime}
||_{C^3(\T^p)}=
||
L_{(\psi,A_\psi,-1)}v_\psi-L_{(\psi^\prime,A_{\psi^\prime},-1)}v_{\psi^\prime}
||_{C^3(\T^p)}\le$$
$$||
L_{(\psi,A_\psi,-1)}(v_\psi-v_{\psi^\prime})
||_{C^3(\T^p)}+
||
(L_{(\psi,A_\psi,-1)}-L_{(\psi^\prime,A_{\psi^\prime},-1)})v_{\psi^\prime}
||_{C^3(\T^p)}\le$$
$$C_9
||
v_\psi-v_{\psi^\prime}
||_{C^0(\T^p)}+
C_{10}(
||
\psi-{\psi^\prime}
||_{C^0(\T^p)}+|A_\psi-A_{\psi^\prime}|
)\cdot
||
v_{\psi^\prime}
||_{C^0(\T^p)}$$
where the the equality comes from (1.15) and the last inequality comes from (1.13) and (1.16). Since the map 
$\fun{}{\psi}{v_\psi}$ is continuous from $Den$ to the $C^0$ topology by step 3, and $\fun{}{\psi}{A_\psi}$ is continuous too (because $\fun{}{\psi}{B_\psi}$ is continuous by step 4 and point 1), $iii$) of this lemma holds), point 2) follows.

\fin

In the next lemma, we show how the fixed points of 
$L_{(\psi,0,-1)}$ induce solutions of $(TS)_{\psi,per}$.

\lem{1.4} 1) (Existence) Given $\psi\in Den$, we can find $A\in\R$ and $\hat v\in C(\T,C^3(\T^p))\cap C^1(\T,C^1(\T^p))$ such that 
$\hat v>0$ and $(\hat v,A)$ solves $(TS)_{\psi,per}$.

\noindent 2) (Uniqueness) Let us suppose that $(\hat v,A)$ and 
$(\hat v^1,A^1)$ are two solutions of $(TS)_{\psi,per}$ with 
$\hat v>0$ and $\hat v^1>0$. Then, $A=A_1$ and 
$\hat v=\l\hat v^1$ for some $\l>0$.

\noindent 3) (Estimates) Let us call $(\hat v_\psi,A_\psi)$ the solution of $(TS)_{\psi,per}$ such that $\hat v_\psi>0$ and 
$\hat v_\psi(0,\cdot)$ satisfies (1.7). Then, there is a constant 
$C_{10}>0$, independent on $\psi\in Den$, such that
$$|A_\psi|+
||
\hat v_\psi
||_{C(\T,C^3(\T^p))}    +
||
\hat v_\psi
||_{C^1(\T,C^1(\T^p))}
\le C_{10}   \eqno (1.17)$$
and
$$\frac{1}{C_{10}}\le\hat v_\psi(t,x)\le C_{10}\qquad
\forall (t,x)\in\T\times\T^p  .   \eqno (1.18)$$

\noindent 4) (Continuous dependence) Let $(\hat v_\psi,A_\psi)$ be as in point 3), and let us consider the map 
$\fun{I}{\psi}{(\hat v_\psi,A_\psi)}$. Then, $I$ is continuous from $Den$ to $[C(\T,C^3(\T^p))\cap C^1(\T,C^1(\T^p))]\times\R$.

\proof As in the proof of lemma 1.3, we set 
$A_\psi=\frac{1}{\b}\log B_\psi$.
For $t\le 0$, we set
$$\hat v_\psi(t,x)=(L_{(\psi,A_\psi,t)}v_\psi)(x)  .  \eqno (1.19)$$
By (1.15), we get that $\hat v_\psi(-1,x)=\hat v(0,x)$; in other words, $\hat v_\psi$ quotients on $\T\times\T^p$; equivalently, it satisfies the second formula of $(TS)_{\psi,per}$. 

Let us prove (1.18) for the function $\hat v_\psi$ defined by (1.19); we prove the inequality on the left, since the one on the right is analogous.

The first equality below is (1.19). Since $\hat v_\psi$ is periodic, we can suppose that $t\in[-1,0]$; now (1.5), implies the first inequality below; the second inequality follows from lemma 1.1 and the fact that 
$t\in[-1,0]$; the third one comes from point 1), $ii$) and $iii$) of lemma 1.3.
$$\hat v_\psi(t,x)=
(L_{(\psi,A_\psi,t)}v_\psi)(x)\ge 
e^{-A_\psi}
e^{\min (t\b P_\psi)}          (\min v_\psi)
E_w\left(
e^{\b\inn{c}{\sqrt\b w(0)}}
\right)  =  $$
$$e^{-A_\psi}
e^{\min (t\b P_\psi)}      (\min v_\psi)
\frac{1}{\sqrt{(2\pi |t|)^p}}
\int_{\R^p}e^{-\frac{|x|^2}{2|t|}}e^{\sqrt{\b}\inn{c}{x}}\dx\ge$$
$$e^{-A_\psi}C_9\min v_\psi\ge\frac{C_9}{C_8}  .  $$
This yields the inequality on the left of (1.18).

We prove (1.17). We begin to note that the estimate on $A_\psi$ follows by point 1), $iii$) of lemma 1.3, and by the fact that 
$A_\psi=\frac{1}{\b}\log B_\psi$. 

We end the proof of (1.17) with the estimates on the derivatives. Let $\tilde w$ be the Brownian bridge with 
$\tilde w(-1)=0=\tilde w(0)$ and let 
$\tilde E_{(0,0)}$ denote its expectation; for $t<0$, let $w$ be the Brownian bridge with $w(t)=0=w(0)$ and let $E_{(0,0)}$ denote its expectation; we recall that
$$w(s)=\frac{1}{\sqrt{|t|}}\tilde w\left(
\frac{s}{|t|}
\right)   .  $$
This yields the second inequality below, while (1.19) yields the first one; the third one comes from the change of variables 
$s=\frac{\tau}{|t|}$. 
$$\hat v_\psi(t,x)=
\left(
\frac{\b}{2\pi|t|}
\right)^\frac{p}{2}
e^{-\b\inn{c}{x}}\cdot$$
$$\int_{\R^p}
e^{
-\frac{|z|^2}{2|t|}
}
e^{
\b\inn{c}{\frac{1}{\sqrt\b}z+x}
}
v_\psi\left(
\frac{1}{\sqrt\b}z+x
\right)
E_{(0,0)}\left[
e^{
\int_{t}^0\b P_\psi(\tau,x+\frac{1}{\sqrt\b}w(\tau)+\frac{\tau+t}{|t|}z)
\dr\tau
}
\right]   \dr z=$$
$$\left(
\frac{\b}{2\pi|t|}
\right)^\frac{p}{2}
e^{-\b\inn{c}{x}}\cdot$$
$$\int_{\R^p}
e^{
-\frac{|z|^2}{2|t|}
}
e^{
\b\inn{c}{\frac{1}{\sqrt\b}z+x}
}
v_\psi\left(
\frac{1}{\sqrt\b}z+x
\right)
\tilde E_{(0,0)}\left[
e^{
\int_{t}^0\b P_\psi(\tau,x+\frac{1}{\sqrt{\b |t|}}w(\frac{\tau}{|t|})+\frac{\tau+t}{|t|}z)
\dr\tau
}
\right]   \dr z=$$
$$\left(
\frac{\b}{2\pi|t|}
\right)^\frac{p}{2}
e^{-\b\inn{c}{x}}\cdot$$
$$\int_{\R^p}
e^{
-\frac{|z|^2}{2|t|}
}
e^{
\b\inn{c}{\frac{1}{\sqrt\b}z+x}
}
v_\psi\left(
\frac{1}{\sqrt\b}z+x
\right)
\tilde E_{(0,0)}\left[
e^{
\int_{-1}^0\b P_\psi(|t|s,x+\frac{1}{\sqrt{\b |t|}}w(s)+(s+1)z)
\dr s
}
\right]   \dr z   .     \eqno (1.20)$$
By point 1), $i$) of lemma 1.3, we can differentiate under the integral sign and get that
$$||\hat v_\psi||_{C([-2,-1],C^3(\T^p))}+
||\hat v_\psi||_{C^1([-2,-1],C^1(\T^p))}   \le C_{10}  .  $$
Since $\hat v_\psi$ is periodic in time, (1.17) follows.

By theorem 9.1 and proposition 6.6 of [6], the Feynman-Kac formula holds for the unbounded final condition $e^{\b\inn{c}{x}}v_\psi$; this, (1.19) and (1.5) imply that $\hat v_\psi$ satisfies the first formula of 
$(TS)_{\psi,per}$ for $t<0$; since it is periodic in $t$, it satisfies it for all times. Moreover, 
$\hat v_\psi>0$ because, by (1.19) and (1.5), it is an integral, with a positive weight, of the positive $v_\psi$. This ends the proof of point 1).

We have just seen that (1.19) gives a bijection between the periodic, positive solutions of $(TS)_{\psi,per}$ and the positive eigenfunctions of $L_{(\psi,0,-1)}$; since the latter are unique up to a multiplicative constant by point 2) of lemma 1.2, we get that the former too are unique up to a multiplicative constant; this proves point 2).

We prove point 4). To prove that the map $\fun{}{\psi}{A_\psi}$ is continuous, it suffices to note that $A_\psi=\frac{1}{\b}\log B_\psi$, that the map $\fun{}{\psi}{B_\psi}$ is continuous  by point 2) of lemma 1.3, and that $B_\psi$ is bounded away from zero and infinity by point 1), $iii$) of the same lemma. 

By point 2) of lemma 1.3, we know that $\fun{}{\psi}{v_\psi}$ is continuous from $Den$ to $C^3(\T^p)$; this and (1.19) easily imply that $\fun{}{\psi}{\hat v_\psi}$ is continuous from $Den$ to 
$C(\T,C^3(\T^p))\cap C^1(\T,C^1(\T^p))$.
 
 \fin
 
 \lem{1.5} 1) (Existence and uniqueness) There is a unique couple $\hat H_\psi(c)\in\R$ and
 $u_\psi\in C(\T,C^3(\T^p))\cap C^1(\T,C^1(\T^p))$ which solves 
 $(HJ)_{\psi,per}$ and satisfies
 $$\int_{\T^p}u_\psi(0,x)\dx=0  .  \eqno (1.21)$$
 
 \noindent 2) (Estimates) There is $C_{11}>0$, independent on 
 $\psi\in Den$, such that, if $u_\psi$ is as in point 1), then
 $$|\bar H_\psi(c)|+
 ||u_\psi||_{C(\T,C^3(\T^p))}+
 ||u_\psi||_{C^1(\T,C^1(\T^p))}\le C_{11}  .  \eqno (1.22)$$
 
 \noindent 3) (Continuous dependence) The couple 
 $(u_\psi,\bar H_\psi(c))$ depends continuously on $\psi$: if 
 $\psi_n\tends\psi$ in $Den$, then 
 $\bar H_{\psi_n}(c)\tends\bar H_\psi(c)$ in $\R$ and 
 $u_{\psi_n}\tends u_\psi$ in 
 $C(\T,C^3(\T^p))\cap C^1(\T,C^1(\T^p))$.
 
 \proof By lemma 1.4, there is a unique couple
 $$(\hat v_\psi,A_\psi)\in
 [C(\T,C^3(\T^p))\cap C^1(\T,C^1(\T^p))]\times\R$$
 which solves $(TS)_{\psi,per}$ and such that 
 $\hat v_\psi(0,\cdot)$ is positive and satisfies (1.7). We have seen at the beginning of this section that, for any $\l>0$, the couple
$$(u_\psi,\bar H_\psi(c))\colon=
(-\frac{1}{\b}\log(\l\hat v_\psi),-\frac{1}{\b}A_\psi)    
\eqno (1.23)$$
solves $(HJ)_{\psi,per}$; vice-versa, if $u_\psi$ solves 
$(HJ)_{\psi,per}$, then its exponential solves $(TS)_{\psi,per}$. Thus, if we define $u_\psi$ as above, for the unique $\l$ for which (1.21) holds, we have existence. Now point 2) of lemma 1.4 implies that all positive solutions of $(TS)_{\psi,per}$ are of the form $(\l\hat v_\psi,A_\psi)$; since we have just seen that there is a bijection between the solutions of $(HJ)_{\psi,per}$ and the positive solutions of $(TS)_{\psi,per}$, we get uniqueness.
 
Formula (1.22) follows from (1.23); indeed, the derivatives of the logarithm of $\hat v_\psi$ are bounded by (1.17) and (1.18). In an analogous way, point 3) follows from point 4) of lemma 1.4.
 
 \fin
 
Let the Lagrangian $\L_{c,\psi}$ be as in the introduction, and let 
 $u_\psi$ be as in lemma 1.5. It is well-known ([9]) that $u_\psi$ satisfies, for $t\le 0$, 
$$u_\psi(t,x)=\min_Y E_w\left\{
\int_t^0\L_{c,\psi}(s,z(s),Y(s,z(s)))\dr s+u_\psi(0,z(0))
\right\}  $$
where $z$ solves the stochastic differential equation
$$\left\{
\eqalign{
{}&\dr z(s)=Y(s,z(s))\dr s+\frac{1}{\sqrt\b}\dr w(s)
\quad s\ge t\cr
{}&z(t)=x
}
\right.     \eqno (SDE)_{t,Y,\d_x}  $$
and $Y(t,z)$ varies among the vector fields continuous in $t$ and Lipschitz in $z$. We have denoted by $E_w$ the expectaction with respect to the Wiener measure. From [9], we get that the minimal $Y_\psi$ is given by 
$$Y_\psi(t,x)=c-\partial_x u_\psi(t,x)  . $$
By (1.22), there is $C_{12}>0$ such that, for any $\psi\in Den$, 
$$||Y_\psi||_{C(\T,C^2(\T^p))}+
||Y_\psi||_{C^1(\T,C(\T^p))}  \le C_{12}  .   \eqno (1.24)$$

\vskip 1pc
\noindent{\bf Definition.} We group in a set $Vect$ all the vector fields $\fun{Y}{\T\times\T^p}{\R^p}$ which satisfy (1.24). The distance on $Vect$ is given by the norm of (1.24).

\vskip 1pc

We would like to consider the law of the  stochastic differential equation above when the initial condition is distributed according to a measure $\mu$. One way to do this is to call $\r_{x_0}$ the solution of $(FP)_{t,Y,\d_{x_0}}$ and to set
$$\r(s,x_0)=\int_{\T^p}\r_{x_0}(s,x)\dr\mu(x_0)  .  $$
Another one, which yields the same law, is to suppose that the Brownian motion is on a probability space $\Omega$ on which there is a random variable $M$ independent on $w(s)$ for 
$s\ge t$ and with law 
$\mu$; we consider the solution $z$ of the stochastic differential equation above with initial condition $M$ and we say that $z$ solves $(SDE)_{t,Y,\mu}$.

\vskip 1pc

Let $Y\in Vect$; by [13], there is $\mu\in C(\T,\Mt)$ which is invariant by the stochastic differential equation; in other words, there is a measure $\mu_0$ such that, if $\mu_t$ is the measure induced by a solution $z$ of $(SDE)_{0,Y,\mu_0}$ for $t\ge 0$, then $\mu_0=\mu_1$. Equivalently, we are saying that there is a weak solution $\mu$ of $(FP)_{Y, per}$. We sketch a proof of this fact: the map which brings the measure $\mu_0$ into 
$\mu_1$, the solution of the Fokker-Planck equation at time $1$, has a fixed point by the Schauder theorem.

We shall use the following classical uniqueness result ([7], proposition 1, [1], theorem 4.1, [17], theorem 5.34) to prove that 
$\mu$ has a smooth density $\rho_Y$.

\lem{1.6} Let $Y\in Vect$. For $i=1,2$, let the map 
$\fun{\nu^i}{[0,+\infty)}{\Mt}$ be continuous and let it be a weak solution of the Fokker-Planck equation, i. e.
$$\int_{\T^p}\phi(0,x)\dr\nu_0^i(x)+
\int_0^{+\infty}\dt\int_{\T^p}\left[
\partial_t\phi+\frac{1}{2\b}\D\phi+\inn{Y}{\partial_x\phi}
\right]   \dr\nu_t^i  =0    \eqno (1.25)$$
for all 
$\phi\in C^1_c([0,+\infty)\times\T^p)\cap C([0,+\infty),C^2(\T^p))$. Let $\nu_0^1=\nu_0^2$. Then, $\nu_t^1=\nu_t^2$ for all $t\ge 0$.

\proof We begin to note that $\mu_t\colon=\nu^2_t-\nu^1_t$ satisfies
$$\int_0^{+\infty}\dt\int_{\T^p}\left[
\partial_t\phi+\frac{1}{2\b}\D\phi+\inn{Y}{\partial_x\phi}
\right]   \dr\mu_s  =0   \qquad
\forall\phi\in C^1_c([0,+\infty)\times\T^p)\cap
C([0,+\infty),C^2(\T^p)).   \eqno (1.26)$$
We have to prove that $\mu_t=0$ for all $t\ge 0$. 

We define the operator $A_Y$ as
$$A_Y\phi=\frac{1}{2\b}\D\phi+\inn{Y}{\partial_x\phi}  .  $$
Let 
$\g\in C^1_c([0,+\infty)\times\T^p)$ and let $t$ be so large that 
$\supp\g\subset\subset[0,t)\times\T^p$. The heat equation with time reversed and final condition in $t$ 
$$\left\{
\eqalign{
{}&\partial_s\phi+A_Y\phi=\g\quad s<t\cr
{}&\phi(t,x)=0\quad\forall x\in\T^p  
}
\right.$$
has a unique solution $\phi$.  We set 
$$\psi(s,x)=\left\{
\matrix{
&\phi(s,x) &s\le t\cr
&0 &s>t     
}    \right.    $$ 
and we see that 
$\psi\in C^1_c([0,+\infty)\times\T^p)\cap C([0,+\infty),C^2(\T^p))$. Indeed, $\psi$ is $C^1$ in $t$ and $C^2$ in $x$ on $s<t$ by theorem 9 of chapter 1 of [10]; it is obviously $C^2$ on $s>t$; it is $C^2$ also in a neighbourhood of $s=t$, because, by the uniqueness of the Cauchy problem for the equation $\partial_s\phi+A_Y\phi=\g$, and the fact that $\g(s,x)=0$ for $s\in[t-\e,t]$ we have that $\phi(s,x)=0$ for $s\ge t-\e$. We use $\psi$ as a test function in (1.26), getting the second equality below.
$$0=\int_0^{+\infty}\dr s\int_{\T^p}[
\partial_s\psi+A_Y\psi
]  \dr\mu_s(x)=
\int_0^{+\infty}\dr s\int_{\T^p}
\g(s,x)\dr\mu_s(x)  .  $$
Since the formula above holds for all 
$\g\in C^1_c([0,+\infty)\times\T^p)$, we get the thesis.

\fin

\lem{1.7} Let $Y\in Vect$. By [13], there is $\mu\in C(\T,\Mt)$ which solves $(FP)_{Y,per}$ in the weak sense. Then, the following holds.

\noindent 1) The measure $\mu$ has density $\r_Y\in Den$.

\noindent 2) The measure $\mu$ is unique.

\noindent 3) There is $C_{13}>0$, independent on $Y\in Vect$, such that
$$||\r_Y||_{C^1(\T\times\T^p)}+
||\r_Y||_{C(\T,C^2(\T^p))}  \le C_{13}  .   $$

\noindent 4) If $Y_n\in Vect$ for all $n$, if $Y\in Vect$ and 
$Y_n\tends Y$ in $C(\T\times\T^p)$, then $\r_{Y_n}\tends \r_Y$ in $C(\T\times\T^p)$.

\proof Classical results about PDE's (see lemma 2.3 below for more details) imply that there is a density $\r_{x_0}$, smooth on 
$(0,+\infty)\times\T^p$, which solves
$$\left\{
\eqalign{
{}&\frac{1}{2\b}\D \r_{x_0}-
{\rm div}[\r_{x_0}\cdot Y]-\partial_t\r_{x_0}=0\cr
{}&\r_{x_0}(t,\cdot)\L^p\tends\d_{x_0} \txt{as}t\tends 0   
}
\right.      \eqno (FP)_{0,Y,\d_{x_0}}$$
where $\L^p$ denotes the Lebesgue measure on $\T^p$.
It is standard that, for $t>0$, $\r_{x_0}(t,\cdot)$ satisfies properties d2) and d3) of the introduction, and that 
$$||\r_{x_0}||_{C^1([1,2]\times\T^p)}+
||\r_{x_0}||_{C([1,2],C^2(\T^p))}\le C_{13}   \eqno (1.27)$$
for a constant $C_{13}>0$ which depends only on the $C^1$ norm of $Y$; as a consequence, $C_{13}$ is the same for all 
$Y\in Vect$ and $x_0\in\T^p$ (again, we refer the reader to lemma 2.3 below).

For $t>0$, we define
$$\r_Y(t,z)=\int_{\T^p}\r_{x_0}(t,z)\dr\mu_0(x_0)  .  
\eqno (1.28)$$
By lemma 1.6, $\r_Y(t,\cdot)$ is the density of $\mu_t$; since 
$\mu_t$ is periodic, we get that 
$\r_Y(t+1,\cdot)=\r_Y(t,\cdot)$.
Point 3) follows from this, (1.27), (1.28) and the fact that norms are convex. One consequence of point 3) is that $\r_Y$ also satisfies hypothesis d1) of the introduction; since we saw above that it satisfies d2) and d3), we get that $\r_Y\in Den$. Again from point 3), we get that $\r_Y$ is a classical solution of $(FP)_{Y,per}$; since by [13] there is only one of them, we get point 2).

We prove point 4). Let $Y_n\tends Y$ in $C(\T\times\T^p)$, and let 
$\r_{Y_n}$ and $\r_Y$ solve $(FP)_{Y_n,per}$ and $(FP)_{Y,per}$ respectively. We have just proved that $\r_{Y_n}$ satisfies point 3) of the thesis; thus, we can apply Ascoli-Arzel\`a\ and get that, up to subsequences, $\r_{Y_n}\tends\r$ in $C(\T\times\T^p)$. Taking limits in (1.25) we see that $\r$ is a weak, periodic solution of 
$(FP)_{Y,per}$; by the uniqueness of point 2), we get that 
$\r=\r_Y$. Thus, any subsequence of $\r_{Y_n}$ has a sub-subsequence converging to $\r_Y$ in $C(\T\times\T^p)$; by a well-known principle, this implies that $\r_{Y_n}\tends\r_Y$ in 
$C(\T\times\T^p)$.

\fin

\vskip 1pc
\noindent{\bf Definition.} Let $C_{13}$ be as in lemma 1.7. 
We group in a set $\r\in Den^{reg}$ the elements of $Den$ which belong to $Lip(\T\times\T^p)$ and such that $||\r||_{Lip(\T\times\T^p)}\le C_{13}$. By point 3) of lemma 1.7, if $Y\in Vect$, then $\r_Y\in Den^{reg}$.

\lem{1.8} There is a continuous map $\fun{\Phi}{Den}{Den}$ whose fixed points $\r_\b$ induce solutions 
$(u_{\r_\b},\r_\b,\bar H_{\r_\b}(c))$ of 
$(HJ)_{\r_\b,per}-(FP)_{c-\partial_x u_{\b}, per}$. Moreover, 
$\Phi(Den)\subset Den^{reg}$.

\proof We define the map $\Phi$ by composition. By lemma 1.5 and formula (1.24), we know that there is a map
$$\fun{\Phi_1}{Den}{Vect,\times\R},\qquad
\fun{\Phi_1}{\psi}{(c-\partial_xu_{\psi},\bar H_\psi(c))}  .  $$
This map is continuous by point 3) of lemma 1.5.

Let $\r_Y$ be as in point 1) of lemma 1.7; by point 4) of this lemma, the map
$$\fun{\Phi_2}{Vect\times\R}{Den},\qquad
\fun{\Phi_2}{(Y,\l)}{\r_Y}$$
is continuous; by point 3), it has image in $Den^{reg}$. Thus, the map $\Phi\colon=\Phi_2\circ\Phi_1$ is continuous from $Den$ to 
$Den$, and has image in $Den^{reg}$, as we wanted.

Let now $\r_\b$ be a fixed point of $\Phi$; we recall that 
$\Phi_1(\r_\b)=(c-\partial_x u_{\r_\b},\bar H_{\r_\b(c)})$, with 
$(u_{\r_\b},\bar H_{\r_\b(c)})$ which satisfies 
$(HJ)_{\r_\b,per}$ and (1.21). Moreover, 
$$\r_\b=\Phi_2\circ\Phi_1(\r_\b)=
\Phi_2(c-\partial_x u_{\r_\b},\bar H_{\r_\b(c)})$$
solves, by the definition of 
$\Phi_2$, $(FP)_{c-\partial_xu_{\r_\b},per}$; in other words, 
$(\r_\b,u_{\r_\b},\bar H_{\r_\b}(c))$ solves 
$(HJ)_{\r_\b,per}-(FP)_{c-\partial_xu_{\r_\b},per}$, and we are done.

\fin

\noindent{\bf Proof of theorem 1.} We begin to show that there are couples $(u_\b,\r_\b)$ which satisfy 
$(HJ)_{\r_\b,per}-(FP)_{c-\partial_x u_{\b}, per}$. By lemma 1.8, this follows if we show that $\Phi$ has fixed points. But this is true by Schauder's fixed point theorem: indeed, by lemma 1.8, $\Phi$ is a continuous map from $Den$ to itself which preserves the compact, convex set $Den^{reg}$.

Let us now call ${\bf S}$ the set of the triples $(u,\r,H)$ such that 
$\r\in Den$ is a weak solution of $(FP)_{c-\partial_x u,per}$ and 
$(u,H)$ is a classical solution of $(HJ)_{\r ,per}$. Let 
$(u^n,\r^n,H^n)\in{\bf S}$ be such that 
$$\int_{\T\times\T^p}\L_{c,\2\r^n}(t,x,\partial_xu^n)\r^n(t,x)\dt\dx
\tends
\inf_{(u,\r,H)\in{\bf S}}
\int_{\T\times\T^p}\L_{c,\2\r}(t,x,\partial_xu)\r(t,x)\dt\dx    .    
  \eqno (1.29)$$
By lemma 1.1,  $\L_{c,\2\r}$ is bounded from below independently on $\r$; as a consequence, the $\inf$ in the right hand side of (1.29) is finite. Note that, if  $\r^n\in Den$, lemma 1.5 implies that 
$c-\partial_xu^n\in Vect$; since $\r^n$ is a fixed point, we get by lemma 1.7 that $\r^n\in Den^{reg}_m$; since $Den^{reg}$ is compact in $Den$, we can suppose that, up to subsequences, 
$$\r^n\tends\bar\r\txt{in}Den   .  $$
By point 3) of lemma 1.5, this implies that 
$$(u^n,H^n)\tends(\bar u,\bar H) \txt{in} 
[C(\T,C^3(\T^p))\cap C^1(\T,C^1(\T^p))]\times\R  ,  $$ 
with $(\bar u,\bar H)$ solving $(HJ)_{\bar\r,per}$. This and point 4) of lemma 1.7 yield that $\bar\r=\r_{c-\partial_x\bar u}$ solves 
$(FP)_{c-\partial_x\bar u,per}$ and satisfies the estimate of point 3) of that lemma.
In other words, $(\bar u,\bar\r,\bar H)\in{\bf S}$; now (1.29) and the last three formulas easily imply that
$$\int_{\T\times\T^p}
\L_{c,\2\bar\r}(t,x,\partial_x\bar u)\bar\r(t,x)\dt\dx=
\inf_{(u,\r,H)\in S}
\int_{\T\times\T^p}\L_{c,\2\r}(t,x,\partial_xu)\r(t,x)\dt\dx  $$
yielding the thesis.

\fin

\vskip 2pc

\centerline{\bf \S 2}
\centerline{\bf The evolution equation}
\vskip 1pc

In this section, we shall prove theorems 2 and 3. We begin with some notation.

We recall that the map
$$\fun{}{(\mu,\nu)}{d_1(\mu,\nu)}$$
is convex, i. e.
$$d_1((1-\l)\nu_1+\l\mu_1,(1-\l)\nu_2+\l\mu_2)\le
(1-\l)d_1(\nu_1,\mu_1)+\l d_1(\nu_2,\mu_2)   .  $$
Indeed, the dual formulation 
$$d_1(\mu,\nu)=\sup\left\{
\int_{\T^p}f\dr(\mu-\nu)\st f\in Lip^1(\T^p)
\right\}    $$
implies that $d_1$ is the supremum of a family of linear functions. 
Since the functions $f$ in the dual formulation belong to  
$Lip_1(\T^p)$ and $\T^p$ has diameter $\sqrt p$, we can as well  suppose that $||f||_\infty\le\2\sqrt p$; as a consequence,
$$d_1(\mu,\nu)\le\sqrt{p}||\mu-\nu||_{tot}  ,  \eqno (2.1)$$
where $||\cdot||_{tot}$ denotes total variation.

\vskip 1pc

\noindent{\bf Definition.} We are going to denote by the norm symbol the distance on $C([-m,0],\Mt)$, which is no norm at all: if 
$R_1,R_2\in C([-m,0],{\cal M}_1(\T^p))$, then we set
$$||
R_1-R_2
||_{C([-m,0],{\cal M}_1(\T^p))}=
\sup_{t\in[-m,0]}d_1(R_1(t),R_2(t))  .  $$
Though this is no norm, it is convex thanks to the convexity of 
$d_1$:
$$||
(1-\l)R_1+\l R_2-(1-\l)\tilde R_1-\l\tilde R_2
||_{C([-m,0],{\cal M}_1(\T^p))}\le$$
$$(1-\l)||
R_1-\tilde R_1
||_{C([-m,0],{\cal M}_1(\T^p))}  +
\l ||
R_2-\tilde R_2
||_{C([-m,0],{\cal M}_1(\T^p))}   .   \eqno (2.2)$$

\vskip 1pc

\vskip 1pc
\noindent{\bf Definition.} For $\mu\in\Mt$ and $m\in\N$, we group in a set $Den_m(\mu)$ all the maps $R\in C([-m,0],\Mt)$ such that 
$R(-m)=\mu$. This space inherits the distance of $C([-m,0],\Mt)$.
\vskip 1pc

\lem{2.1} Let $f\in C^3(\T^p)$ and  let 
$H^Z(t,q,p)=\2|p|^2+Z(t,q)$, with $Z\in C([-m,0],C^3(\T^p))$. 

\noindent 1) Then, there is a unique solution $u^Z$ of
$$\left\{
\eqalign{
{}&\frac{1}{2\b}\D u^Z+\partial_tu^Z-
H^Z(t,x,c-\partial_xu^Z)=0,\quad t\in[-m,0]\cr
{}&u^Z(0,x)=f\quad\forall x\in\T^p   .
}
\right.   \eqno (HJ)^{Z}  $$

\noindent 2) There is $C_{13}>0$, only depending on 
$||f||_{C^3(\T^p)}$, 
$||Z||_{C([-m,0],C^3(\T^p))}$ and $m$, such that
$$||u^Z||_{C^1([-m,0],C^1(\T^p))}+
||u^Z||_{C([-m,0],C^3(\T^p))}\le C_{13}  .   $$

\noindent 3) The map
$$\fun{}{Z}{u^Z}$$
is continuous from $C([-m,0],C^3(\T^p))$ to 
$C([-m,0],C^3(\T^p))\cap C^1([-m,0],C^1(\T^p))$.

\proof We know that the twisted Schroedinger equation with potential $Z$ and final condition $e^{-\b f}\in C^3(\T^p)$ has a unique solution $v^Z$, which can be represented by the Feynman-Kac formula (1.20) with $e^{-\b f}$ in stead of $v_\psi$ and $Z$ in stead of $P_\psi$. Since $e^{-\b f}>0$, we get that $v^Z>0$ too. We saw in section 1 that $u^Z=-\frac{1}{\b}\log v^Z$ solves 
$(HJ)^Z$, and point 1) follows. 

Points 2) and 3) follow as in section 1 if we prove that
$$||\partial_x v^Z||_{C^1([-m,0],C^1(\T^p))}+
||v^Z||_{C([-m,0],C^3(\T^p))}\le C_{14} $$
and that the map $\fun{}{Z}{(v^z,\partial_x v^z)}$ is continuous. Since $e^{-\b f}$, the final condition of the Schroedinger equation, is of class $C^3$, this is a standard result; for instance, differentiation under the integral sign in (1.20) gives the estimate on 
$||v^Z||_{C([-m,0],C^3(\T^p))}$; from this and the fact that $v^Z$ solves the Scroedinger equation, we get the estimate on 
$||v^Z||_{C^1([-m,0],C^1(\T^p))}$.

\fin

Recalling lemma 1.1, we get this immediate consequence.

\cor{2.2} 1) Let $f\in C^3(\T^p)$, let $\m\in\Mt$ and let 
$R\in Den^m(\m)$. Then, there is a unique solution $u_R$ of
$$\left\{
\eqalign{
{}&\frac{1}{2\b}\D u_R+\partial_tu_R-
H_R(t,x,c-\partial_xu_R)=0,\quad t\in[-m,0]\cr
{}&u_R(0,x)=f\quad\forall x\in\T^p   .
}
\right.   \eqno (HJ)_{R, f}  $$

\noindent 2) There is $C_{14}=C_{14}(m)>0$, independent of 
$\m\in\Mt$ and on $R\in Den_m(\m)$, such that
$$||u_R||_{C^1([-m,0],C^1(\T^p))}+
||u_R||_{C([-m,0],C^3(\T^p))}\le C_{14}  .   $$

\noindent 3) The map
$$\fun{}{R}{u_R}$$
is continuous from $Den_m(\m)$ to 
$C([-m,0],C^3(\T^p))\cap C^1([-m,0],C^1(\T^p))$.

\rm

\vskip 1pc
\noindent{\bf Definition.} By point 2) of corollary 2.2, there is 
$C_{15}>0$ such that, setting $Y=c-\partial_xu_R$, we have
$$||Y||_{C^1([-m,0],C(\T^p))}+
||Y||_{C([-m,0],C^2(\T^p))}\le C_{15}$$
with $C_{15}$ independent on $R\in C([-m,0],\Mt)$. We group in a set $Vect_m$ all the vector fields $Y$ on $[-m,0]\times\T^p$ which satisfy the estimate above. The distance on $Vect_m$ is the one induced by the norm above.

\vskip 1pc

\lem{2.3} Let $Y\in Vect_m$, and let $\m\in\Mt$. Then, the following holds.

\noindent1) There is a unique $R_Y\in C([-m,0],\Mt)$ which solves 
$(FP)_{-m,Y,\m}$ in the weak sense.

\noindent 2) For $t\in(-m,0]$, $R_Y(t)$ has density $\r_Y$. There are $\fun{C_{16},C_{17}}{(-m,0]}{[0,+\infty)}$, independent on $Y\in Vect_m$ and on $\m\in\Mt$, such that

\noindent a) $C_{17}(T)\tends 0$ as $T\tends -m$, $C_{16}$ and 
$C_{17}$ are bounded on $(-m+\e,0]$ for all $\e\in(0,m)$ and

\noindent b) For $T\in(-m,0]$, we have
$$\left\{
\eqalign{
{}&||\r_Y||_{C^1((T,0)\times\T^p)}+
||\r_Y||_{C((T,0),C^2(\T^p)}\le C_{16}(T)\cr
{}&d_1(R_Y(T),\m)\le C_{17}(T)     
}
\right.   \eqno (2.3)$$
where $d_1$ denotes the $1$-Wasserstein distance.

\proof The uniqueness of point 1) comes from lemma 1.6; for the existence, we begin to recall from PDE theory (see for instance chapter 1 of [10]) that, for 
$x_0\in\T^p$, $(FP)_{-m,Y,\d_{x_0}}$ has a solution $R_{x_0}$ with density $\r_{x_0}$. Always from [10], the function $\r_{x_0}$ satisfies the first formula of (2.3) for a constant $C_{16}(T)$ which  depends neither on $x_0\in\T^p$ nor on the particular element 
$Y\in Vect_m$. Moreover, as 
$T\tends-m$, we get from [10] that, if $g\in C(\T^p)$, then
$$\int_{\T^p}g(x)\dr R_{x_0}(T)\tends g(x_0)$$
uniformly in $x_0\in\T^p$; since $d_1$ induces the weak$\ast$ topology and $\T^p$ is compact, we have that $d_1(R_{x_0}(T),\d_{x_0})\le C_{17}(T)$, for a constant $C_{17}(T)$ which depends neither on $x_0\in\T^p$ nor on $Y\in Vect_m$, and such that $C_{17}(T)\tends 0$ as 
$T\tends -m$. In other words, $\r_{x_0}$ satisfies (2.3) for two uniform constants $C_{16}(T)$, $C_{17}(T)$, depending neither on $x_0$ nor on $Y\in Vect_m$.

Now we set
$$\r_Y(t,x)=\int_{\T^p}\r_{x_0}(t,x)\dr\m(x_0)  .  \eqno (2.4)$$
Clearly, $\r_Y$ is a solution of $(FP)_{-m,Y,\m}$, and this ends the proof of point 1). 

We have seen that $\r_{x_0}$ satisfies the first formula of (2.3); since norms are convex, (2.4) implies that $\r_Y$ too satisfies this formula. Now $\r_{x_0}$ satisfies 
$d_1(R_{x_0}(T),\d_{x_0})\le C_{17}(T)$, and the map
$$\fun{}{(\mu,\nu)}{d_1(\mu,\nu)}$$
is convex; it follows again by (2.4) that $\r_Y$ too satisfies the second formula of (2.3).

\fin

\vskip 1pc
\noindent{\bf Definition.} We define $Den_m^{reg}(\m)$ as the subset of the elements $R\in Den_m(\m)$ which, for $t\in(-m,0]$, have a density $\r$ with respect to the Lebesgue measure. Moreover, we ask that $R$ and $\r$ satisfy
$$\left\{
\eqalign{
{}&||\r||_{Lip([T,0]\times\T^p)}\le C_{16}(T),\quad
\forall T\in(-m,0]\cr
{}&d_1(R(T),\m)\le C_{17}(T),\quad 
\forall T\in(-m,0]    
}
\right.  \eqno (2.5)$$
where $C_{16}(T)$ and $C_{17}(T)$ are the same two constants of (2.3). By lemma 2.3, if $Y\in Vect$, $\mu\in{\cal M}_1(\T^p)$ and $R_Y$ solves $(FP)_{-m,Y,\mu}$ in the weak sense, then 
$R_Y\in Den^{reg}_m(\mu)$.

\lem{2.4} $Den^{reg}_m(\m)$ is compact in $Den_m(\m)$ for the 
$C([-m,0],\Mt)$ topology.

\proof Let $R_n\in Den^{reg}_m(\m)$ have density $\r_n$ for 
$n\in\N$. We must show that it has a subsequence converging in 
$Den_m(\mu)$.

Since $\r_n$ satisfies the first formula of (2.5), Ascoli-Arzel\`a\ implies that, up to subsequences, 
$\r_n\tends\r$ in $C^0_{loc}((-m,0]\times\T^p)$; clearly, 
$\r$ satisfies the first formula of (2.5). Denoting by $\L^p$ the Lebesgue measure on $\T^p$, we set $R(t)=\r(t)\L^p$ and we see that, for any fixed $T\in(-m,0]$,
$$d_1(
R_n(T),R(T)
)\le
\sqrt{p}||R_n(T)-R(T)||_{tot}=
\sqrt{p}||\r_n(T)-\r(T)||_{L^1(\T^p)}
\tends 0\txt{as}n\tends+\infty     $$
where the first inequality comes from (2.1) and the limit from the fact that $\r_n\tends\r$ in $C^0_{loc}((-m,0]\times\T^p)$. 
Since $R_n$ satisfies the second formula of (2.5), we have that
$$d_1(
R_n(T),\m
)   \le 
C_{17}(T),\quad\forall T\in(-m,0],\quad\forall n\ge 1  .  $$
The last two formulas imply that $R$ satisfies the second formula of (2.5).

It remains to prove that $R_n\tends R$ in 
$C([-m,0],{\cal M}_1(\T^p))$; it suffices to note that, for 
$\d\in(0,m)$, 
$$\sup_{t\in[-m,0]}   d_1(R_n(t),R(t))  \le$$
$$\sup_{t\in[-m,-m+\d]}[
d_1(R_n(t),\mu)+d_1(\mu,R(t))
]+
\sup_{t\in[-m+\d,0]}d_1(R_n(t),R(t))\le$$
$$2C_{17}(-m+\d)+
\sqrt{p}\sup_{t\in[-m+\d,0]}||\r_n(t)-\r(t)||_{L^1(\T^p)}    $$
where the last inequality comes from the second formula of (2.5) and from (2.1). Since $C_{17}(T)\tends 0$ as $T\searrow-m$, we can fix $\d>0$ so that the first term on the right is smaller than 
$\e$; having thus fixed $\d$, we take $n$ so large that, by convergence in 
$C^0_{loc}((-m,0]\times\T^p)$, the second term on the right is smaller than $\e$, and we are done.

\fin

We only sketch the proof of the next lemma, since it is identical to point 4) of lemma 1.7.

\lem{2.5} Given $\e>0$, we can find $\d>0$ with the following property. Let $Y,\bar Y\in Vect_m$ and let 
$\mu\in{\cal M}_1(\T^p)$; let $R_{\bar Y}$ and 
$R_Y$ satisfy $(FP)_{-m,\bar Y,\mu}$ and $(FP)_{-m,Y,\mu}$ respectively. Let $||\bar Y- Y||_{C([-m,0]\times\T^p)}\le\d$. Then, 
$||R_{\bar Y}- R_Y||_{C([-m,0],\Mt)}\le\e$.

\proof Let $\{ Y_n \}_{n\ge 1}$, $\{ \bar Y_n \}_{n\ge 1}$ be two sequences in $Vect_m$ and let $\{ \mu_n \}_{n\ge 1}\subset\Mt$. We suppose that 
$||\bar Y_n- Y_n||_{C([-m,0]\times\T^p)}\tends 0$; we let  
$R_{Y_{n}}$ solve $(FP)_{-m,Y_n,\mu_n}$ and  
$R_{\bar Y_{n}}$ solve $(FP)_{-m,\bar Y_n,\mu_n}$; we have to prove that
$$||R_{\bar Y_n}- R_{Y_n}||_{C([-m,0],\Mt)}\tends 0   .  $$
Let us suppose by contradiction that this does not hold; then there is $\e>0$ and a subsequence (which we denote by the same index) such that 
$$||R_{\bar Y_n}- R_{Y_n}||_{C([-m,0],\Mt)}>\e \qquad\forall n  .  
$$
Since $\bar Y_n,Y_n\in Vect_m$ and 
$||\bar Y_n- Y_n||_{C([-m,0]\times\T^p)}\tends 0$, by 
Ascoli-Arzel\`a\ up to taking subsequences we can suppose that 
$\bar Y_n,Y_n\tends Y$ in 
$C([-m,0]\times\T^p)$; we can also suppose that 
$\mu_n\tends\mu$. To reach a contradiction with the formula above, it suffices to show that $R_{Y_n}$ and $R_{\bar Y_n}$ both converge to $R_Y$; since the proof for $R_{\bar Y_n}$ is analogous, we prove convergence for $R_{Y_n}$.

We note that $\{ R_{Y_{n}} \}$ in contained in 
$Den_m^{reg}(\mu_n)$ by lemma 2.3; thus, by lemma 2.4, it has a subsequence converging to a limit $R$. Since $R_{Y_{n}}$ is a weak solution of $(FP)_{-m,Y_{n},\mu_n}$, we easily get that $R$ is a weak solution of $(FP)_{-m,Y,\mu}$; by lemma 1.6, $R=R_Y$. In other words, every subsequence of $R_{Y_n}$ has a sub-subsequence converging to $R_Y$; this implies that $R_{Y_n}$ converges to $R_Y$, and we are done.

\fin

\noindent{\bf Proof of theorem 2.} For $Q\in Den_m(\mu)$, let 
$u_Q$ be as in corollary 2.2; for $Y\in Vect_m$, let 
$R_Y=\r_Y\L^p$ be as in lemma 2.3. The two maps 
$$\fun{}{Den_m(\mu)}{Vect_m},\qquad
\fun{}{Q}{c-\partial_xu_{Q}}$$ 
and 
$$\fun{}{Vect_m}{Den_m(\mu)},\qquad
\fun{}{Y}{R_{Y}}$$ 
are both continuous: the first one, by point 3) of of corollary 2.2, the second one by lemma 2.5. Let us call $\Phi$ their composition:
$$\fun{\Phi}{Den_m(\m)}{Den_m(\m)},\qquad
\fun{\Phi}{Q}{R_{c-\partial_xu_Q}}  .  $$
Being the composition of two continuous functions, $\Phi$ is continuous; moreover, by point 2) of lemma 2.3, it has image in 
$Den^{reg}_m(\m)$; this latter set is clearly convex, and it is compact in $Den_m(\m)$ by lemma 2.4. Thus, we have that
$$\fun{\Phi}{Den^{reg}_m(\m)}{Den^{reg}_m(\m)}  .  $$
We apply the Schauder fixed point theorem and we get that $\Phi$ has a fixed point in $Den^{reg}_m(\m)$. With the same argument as in the proof of theorem 1, we see that, if $R$ is a fixed point of 
$\Phi$, then $(u_R,R)$ solves 
$(HJ)_{R,f}-(FP)_{-m,c-\partial_xu_R,\m}$. This yields existence. 

We continue as in the proof of theorem 1. Let us call ${\bf S}$ the set of the couples $(u,R)$ where $u$ is a classical solution of 
$(HJ)_{R,f}$ and $R\in Den_m(\mu)$ is a weak solution of 
$(FP)_{-m,c-\partial_x u,\mu}$. 

Let us consider a sequence 
$(u_n,R_n)\in{\bf S}$ such that, denoting by $\r_n$ the density of $R_n$,
$$\int_{-m}^0\dr t\int_{\T^p}\L_{c,\2 R_n}(t,x,c-\partial_x u_n)\r_n\dr t\dx \tends
\inf_{(u,R)\in{\bf S}} 
\int_{-m}^0\dr t\int_{\T^p}\L_{c,\2 R}(t,x,c-\partial_x u)\r\dr t\dx  . $$
Whatever is $R_n\in Den_m(\mu)$, $u_n$ satisfies the estimates of point 2) of corollary 2.2; in particular, 
$c-\partial_x u_n\in Vect_m$. Since $R_n$ satisfies 
$(FP)_{-m,c-\partial_xu_n,\mu}$ lemma 2.3 implies that 
$R_n\in Den_m^{reg}(\mu)$; by lemma 2.4, up to subsequences we can suppose that 
$R_n\tends\bar R$, with $\bar R\in Den^{reg}_m(\mu)$. By point 3) of corollary 2.2, we get that $u_n\tends\bar u$ in 
$C^1([-m,0],C^1(\T^p))\cap C([-m,0],C^3(\T^p))$, and that $\bar u$ solves 
$(HJ)_{\bar R,f}$. Thus, $(\bar u,\bar R)\in{\bf S}$; now, the formula above easily implies that $(\bar u,\bar R)$ is minimal in 
${\bf S}$.

\fin

We turn to the proof of theorem 3; our route will pass through an approximation with a finite number of particles.

\noindent{\bf Definitions.} Let us define the Lagrangian for one particle as
$$\fun{L_c}{\T\times\T^p\times\R^p}{\R},\qquad
L_c(t,x,y)=\2|y|^2-\inn{c}{y}-V(t,x)  .  $$
The Lagrangian for $n$ particles, each of mass $\frac{1}{n}$, is 
$$\fun{L^n_c}{\T\times(\T^p)^n\times(\R^p)^n}{\R}$$
$$L^n_c(t,(x_1,\dots,x_n),(y_1,\dots,y_n))=
\frac{1}{n}\sum_{i=1}^nL_c(t,x_i,y_i)+
\frac{1}{2n^2}\sum_{i,j=1}^nW(x_i-x_j)  .  $$
Let $U$ be as in the statement of theorem 2.
For any given $z=(z_1,\dots,z_n)\in(\T^p)^n$, we define
$$U^n(-m,z)=$$
$$\inf E_{w_1,\dots,w_n}\left\{
\int_{-m}^0L^n_c(s,X^n(-m,s,z),Y^n(s,X^n(-m,s,z)))\dr s
+U(R^n(-m,0))
\right\}    .  \eqno (2.6)$$
The infimum above is over all vector fields 
$Y^n(s,x)=(Y^n_1(s,x_1),Y^n_2(s,x_2),\dots,Y^n_n(s,x_n))$ continuous in $s$ and Lipschitz in $x$; each component of the function 
$$X^n(-m,s,z)=
(X^n_1(-m,s,z_1),X^n_2(-m,s,z_2),\dots,
X^n_n(-m,s,z_n))      \in(\R^p)^n$$ 
solves the stochastic differential equation on $\R^p$
$$\left\{
\eqalign{
{}&\dr X^n_i(-m,s,z_i)=Y^n_i(s,X^n_i(-m,s,z_i))\dr t+\dr w_i(s)\quad 
s\ge -m,\quad i\in(1,\dots,n)\cr
{}&X^n_i(-m,-m,x_i)=z_i  .  
}
\right.     \eqno (SDE)_{t,Y^n_i,\d_{z_i}}$$
In the formula above, each $w_i$ is a standard Brownian motion on $\R^p$; the $w_i$ are independent and $E_{w_1,\dots,w_n}$ denotes the expectaction with respect to the product of the Wiener measures. It remains to define $R^n(-m,0)$; to do this, we let 
$\r_i^n(-m,s,x)$  be the density on $\T^p$ which solves 
$(FP)_{-m,Y_{i}^n,\d_{z_i}}$ and we set, for $t\in[-m,0]$,
$$\r^n(-m,t,x)=\frac{1}{n}\sum_{i=1}^n\r_i^n(-m,t,x),\qquad
R^n(-m,t)=\r^n(-m,t,x)\L^p  .   \eqno (2.7)$$
We note that we are not considering the most general vector field $Y$ on $(\T^p)^n$. On the contrary, we assign to each particle $x_i\in\T^p$ a control $Y_i$ which depends only on $x_i$, and not on the positions of the other particles; these, however, interact with $x_i$ via the potential $W$. We have chosen this particular problem because we want $U_n(-m,z)$ to converge, as $n\tends+\infty$, to 
$\Lambda^m_cU$; we recall that, in the definition of 
$\Lambda_c^m$, there is a control $Y$ which depends on the single particle in $\T^p$.

\lem{2.6} Let us suppose that $U$ is as in the statement of theorem 1 and let $U^n(-m,z)$ be defined as in (2.6). Then for any fixed $n\in\N$, the infimum in (2.6) is a minimum.

\proof Let $\{ Y^{n,k} \}_{k\ge 1}$ be a minimizing sequence. We are going to show that we can build another minimizing sequence, say $\{ \tilde Y^{n,k} \}_{k\ge 1}$, which is  Lipschitz in $(t,x)$ uniformly in $k$. Once we know this, the lemma follows by Ascoli-Arzel\`a.

For the vector field $Y^{n,k}$, let us define $\r_i^{n,k}$ and 
$\r^{n,k}$ as in (2.7); we set
$$\fun{L^n_{c,Y^{n,k},i}}{[-m,0]\times\T^p\times\R^p}{\R}  ,  $$
$$L^n_{c,Y^{n,k},i}(s,x,\dot x)=
L_c(s,x,\dot x)-
\frac{1}{n}\sum_{j\not=i}
\int_{\T^p}W(x-y)\r_j^{n,k}(-m,s,y)\dr y  .  \eqno (2.8)$$
Note that, in contrast with $L^n_c$, a factor $\2$ in the interaction sum is missing. 
We know from lemma 1.1 that the potential in $L^n_{c,Y^{n,k},i}$ satisfies a uniform $C^3$ estimate. By [9], for 
$(t,x)\in[-m,0]\times\T^p$, there is 
$\tilde Y_i^{n,k}$ on which the minimum below is attained
$$u^{n,k}_i(t,x)\colon=
\min E_w\left\{
\int_t^0L^n_{c,Y^{n,k},i}(s, X, Y)\dr s
+f(X(t,0,x))
\right\}     \eqno (2.9)$$
with $X(t,s,x)$ which solves $(SDE)_{t,Y,\d_{x}}$; the minimum is taken over all the Lipschitz vector fields $Y$. Always by [9], 
$\tilde Y_i^{n,k}=c-\partial_xu^{n,k}_i(t,x)$ and $u^{n,k}_i$ solves the Hamilton-Jacobi equation for the Lagrangian 
$L^n_{c,Y^{n,k},i}$ and final condition $f$. By lemma 2.1, 
$$||u_i^{n,k}||_{C^1([-m,0],C^1(\T^p))}+
||u^{n,k}_i||_{C([-m,0],C^3(\T^p))}   $$
is bounded in terms of the $C^3$ norm of the potential of 
$L^n_{c,Y^{n,k},i}$. By lemma 1.1, the latter depends neither on $n$ nor on $k$; thus, $\tilde Y^{n,k}_i$ belongs to 
$Vect_m$; in particular, it is Lipschitz uniformly in $n$ and $k$.

In the following, whenever we have a drift, say $Y_i^B$, we shall denote by $X_i^B(t,s,x_i)$ the solution of 
$(SDE)_{t,Y_i^B,\d_{x_i}}$; we shall set 
$X^B=(X_1^B,\dots,X_n^B)$ and $z=(z_1,\dots,z_n)$.

We are going to isolate the first particle and show that the mean action decreases if we substitute $Y^{n,k}_1$ with the smoother 
$\tilde Y^{n,k}_1$ defined above. Since the interaction potential is even and satisfies $W(0)=0$, we get the first equality below; since the Brownian motions $(w_1,\dots,w_n)$ are independent, we get the second one.
$$E_{w_1,\dots,w_n}\left\{
\int_{-m}^0L^n_c(s,X^{n,k}(-m,s,z),Y^{n,k}(s,X^{n,k}(-m,s,z)))\dr s+
U(R^n(-m,0))
\right\}  =$$
$$E_{w_1,\dots,w_n}\bigg\{
\frac{1}{n}\sum_{j\not =1}
\int_{-m}^0
L_c(s,X^{n,k}_j(-m,s,z_j),Y^{n,k}_j(s,X^{n,k}(-1,s,z_j)))\dr s+$$
$$\frac{1}{2n^2}\sum_{j,i\not=1}
\int_{-m}^0W(X_i^{n,k}(-m,s,z_i)-X_j^{n,k}(-m,s,z_j))\dr s
+\frac{1}{n}\sum_{j\not=1}f(X_j(-m,0,z_j))
\bigg\}  +$$
$$E_{w_1,\dots,w_n}\bigg\{
\frac{1}{n}\int_{-m}^0[
L_c(s,X^{n,k}_1(-m,s,z_1),Y^{n,k}_1(s,X^{n,k}_1(-m,s,z_1)))+$$
$$\frac{1}{n^2}\sum_{j\not=1}
\int_{-m}^0W(X_1^{n,k}(-m,s,z_1)-X_j^{n,k}(-m,s,z_j))
]\dr s
+\frac{1}{n}f(X_1(-m,0,z_1))
\bigg\} = $$
$$E_{w_1,\dots,w_n}\bigg\{
\frac{1}{n}\sum_{j\not =1}
\int_{-m}^0L_c(s,X_j^{n,k}(-m,s,z_j),Y^{n,k}_j(s,X_j^{n,k}(-m,s,z_j)))\dr s+   \leqno a1)$$
$$\frac{1}{2n^2}\sum_{j,i\not=1}
\int^0_{-m}W(X_i^{n,k}(-m,s,z_i)-X_j^{n,k}(-m,s,z_j))\dr s
+\frac{1}{n}\sum_{j\not=1}f(X_j(-m,0,z_j))
\bigg\}  +  \leqno a2)$$
$$\frac{1}{n}E_{w_1}\left\{
\int_{-m}^0L^n_{c,Y^{n,k},1}(s,X_1^{n,k}(-m,s,z_1),
Y^{n,k}_1(s,X_1^{n,k}(-m,s,z_1)))\dr s
+f(X_1(-m,0,z_1))
\right\}   .  \leqno a3)$$
If we consider $(\tilde Y^{n,k}_1,Y^{n,k}_2,\dots,Y^{n,k}_n)$ instead of $(Y^{n,k}_1,Y^{n,k}_2,\dots,Y^{n,k}_n)$, we see that the terms $a1)$ and $a2)$ in the formula above remain the same, while, by our choice of $\tilde Y^{n,k}_1$, $a3)$  gets smaller. After applying this procedure to each coordinate, we get a sequence 
$\tilde Y^{n,k}=(\tilde Y^{n,k}_1,\dots,\tilde Y^{n,k}_n)$ which satisfies the following two properties.
$$E_{w_1,\dots,w_n}\left\{
\int_{-m}^0
L^n_c(s,\tilde X^{n,k}(-m,s,x),
\tilde Y^{n,k}(s,\tilde X^{n,k}(-m,s,x)))\dr s
+\frac{1}{n}\sum_{j=1}^n f(\tilde X^{n,k}(-m,s,x))
\right\} \le        \leqno\bullet  $$
$$E_{w_1,\dots,w_n}\left\{
\int_{-m}^0L^n_c(s,X^{n,k}(-m,s,x),Y^{n,k}(s,X^{n,k}(-m,s,x)))\dr s
+\frac{1}{n}\sum_{j=1}^n f(X^{n,k}(-m,s,x))
\right\}  .  $$

\noindent $\bullet$ $\tilde Y^{n,k}\in (Vect_m)^n$.

In particular, $\tilde Y^{n,k}_i\in Vect_m$ for all $i$, $n$ and $k$; as a consequence, we can apply point 2) of lemma 2.3, getting that $\{ \r_i^{n,k} \}_{n,k}\subset Den^{reg}_m(\d_{z_i})$ for all $i$. By lemma 2.4, we find a subsequence, which we denote by the same index, such that 
$$(\r^{n,k}_1\L^p,\dots,\r^{n,k}_n\L^p)\tends
(\r^{n}_1\L^p,\dots,\r^{n}_n\L^p)  \txt{in} 
Den_m(\d_{z_1})\times\dots\times Den_m(\d_{z_n})  .  $$
Thus, for each $i$,
$$\frac{1}{n}\sum_{j\not=i}
\int_{\T^p}W(x-y)\r^{n,k}_j(-m,s,y)\dr y\tends
\frac{1}{n}\sum_{j\not=i}
\int_{\T^p}W(x-y)\r^{n}_j(-m,s,y)\dr y    $$
in $C([-m,0],C^3(\T^p))$.  By point 3) of lemma 2.1, this implies that
$$(\tilde Y^{n,k}_1,\dots,\tilde Y^{n,k}_n)\tends
(\tilde Y^{n}_1,\dots,\tilde Y^{n}_n)\txt{in}
C^1([-m,0],C(\T^p))\cap C([-m,0],C^2(\T^p))    $$
and that each $\tilde Y^n_i$ is minimal for $L^n_{c,\tilde Y,i}$.
By the last formula and lemma 2.5, we get that $\r^n_i$ solves 
$(FP)_{-m,\tilde Y^n_i,\d_{z_i}}$. The last three formulas imply that  
$$E_{w_1,\dots,w_n}\left\{
\int_{-m}^0L^n_c(s,\tilde X^n(-m,s,z),
\tilde Y^n(s,\tilde X^n(-m,s,z)))\dr s
+U(R^n(-m,0,z))
\right\} =$$
$$\lim_{n\tends+\infty}
E_{w_1,\dots,w_n}\left\{
\int_{-m}^0L^n_c(s,\tilde X^{n,k}(-m,s,z),
\tilde Y^{n,k}(s,\tilde X^{n,k}(-m,s,z)))\dr s
+U(R^{n,k}(-m,0,z))
\right\} .$$
Since
$\{ (\tilde Y^{n,k}_1,\dots,\tilde Y^{n,k}_n) \}_{k\ge 1}$ is a minimizing sequence, we get that 
$(\tilde Y^n_1,\dots,\tilde Y^n_n)$ is minimal.

\fin

From the proof of the last lemma, we extract the following corollary: it says that the minimum in (2.6) is a Nash equilibrium ([4]).  Note one fact about the value function $u^n_i$ in the corollary below: for simplicity, we let $i=1$ . Then the function  $u^n_1$ depends not only on $(x_2,\dots,x_n)$, but on $x_1$ too: namely, if $x_1$ moves, the drifts $(Y_2,\dots,Y_n)$ will adjust, and the Lagrangian 
$L_{c,Y,1}$ will change. If it hadn't been too clumsy, we could have written $u_1^{n,(x_1,\dots,x_n)}$ and said that 
$c-\partial_x u_1^{n,(x_1,\dots,x_n)}(x)$ is the best drift for particle 
$x_1$. 

\cor{2.7} Let 
$\bar Y^n(t,x)=(\bar Y^n_1(t,x_1),\dots,\bar Y^n_n(t,x_n))$ be minimal in (2.6) and let $L^n_{c,\bar Y^n,i}$ be defined as in (2.8).  Let
$$u^n_i(t,x)=\min_Y E_w\left\{
\int_t^0 L^n_{c,Y,i}(s,X,Y)\dr s
+f(X(t,0,x))
\right\}      $$
where $X$ solves $(SDE)_{t,Y,\d_x}$ and the minimum is taken among all Lipschitz vector fields 
$Y$ on $[-m,0]\times\T^p$. Then, for each $i$ we have that 
$\bar Y^n_i(t,x_i)=c-\partial_xu^n_i(t,x_i)$.

\proof If for one $i$ we had $\bar Y^n_i\not=c-\partial_x u^n_i$, then, isolating particle $i$ as in the last lemma, we could see that the vector field 
$$(
\bar Y^n_1,\dots,\bar Y^n_{i-1},c-\partial_xu^n_i,
\bar Y^n_{i+1},\dots,\bar Y^n_n
)$$
has a lower Lagrangian action, contradicting the minimality of 
$\bar Y^n$.

\fin

\lem{2.8} Let $\m\in\Mt$ and let us suppose that
$$\frac{1}{n}(\d_{z_1}+\dots+\d_{z_n})\tends\m
\txt{in}\Mt  .  \eqno (2.10)$$
Let $Y^n=(Y^n_1,\dots,Y^n_n)$ be a drift minimal in (2.6); by corollary 2.7, $Y^n_i=c-\partial_xu^n_i$ for the value function $u^n_i$ defined in (2.9). Let $\r^n$ be defined as in (2.7).

Then, there is $(u,\r)$ which satisfies 
$(HJ)_{\r,f}-(FP)_{-m,c-\partial_xu,\m}$, and a subsequence 
$\{ n_k \}$ such that
$$\left\{
\matrix{
\r^{n_k}\L^p\tends\r\L^p \txt{in}C([-m,0],{\Mt})\cr
{}\cr
\sup_{0\le i\le n_k}
||u^{n_k}_i-\partial_x u||_{C^1([-m,0],C^1(\T^p))}+
||u^{n_k}_i- u||_{C([-m,0],C^3(\T^p))}\tends 0\cr
{}\cr
\sup_{0\le i\le n_k}
||Y^{n_k}_i- (c-\partial_x u)||_{C^1([-m,0],C(\T^p))}+
||Y^{n_k}_i- (c-\partial_x u)||_{C([-m,0],C^2(\T^p))}   \tends 0  .  
}   
\right.     \eqno (2.11)$$
Moreover, the function $U^{n_k}(-m,z_1,\dots,z_{n_k})$ defined in (2.6) converges to the function $U(-m,\m)$ defined by
$$U(-m,\m)\colon =\int_{[-m,0]\times\T^p}
\L_{c,\2\r}(t,x,c-\partial_xu)\r(t,x)\dx\dt+U(\r)  .  \eqno (2.12)$$

\proof Since $Y^n_i=c-\partial_x u^n_i$, the third formula of (2.11) follows from the second one; we prove the first two ones.

\noindent{\bf Step 1.} We prove the convergence of the densities.

For $i,j\in(1,\dots,n)$, we consider the densities
$$\matrix{
\hat\r_i^n(-m,s,x)=\frac{1}{n-1}\sum_{l\not=i}\r_l^n(-m,s,x),\qquad
\hat\r_j^n(-m,s,x)=\frac{1}{n-1}\sum_{l\not=j}\r_l^n(-m,s,x)  \cr
{}\cr
\r^n(-m,s,x)=\frac{1}{n}\sum_{l=1}^n\r^n_l(-m,s,x)
}
\eqno (2.13)$$
where $\r_l^n$ is the same as in formula (2.7). Let 
$R^n_i=\r^n_i\L^p$, $\hat R^n_i=\hat\r^n_i\L^p$ and 
$R^n=\r^n\L^p$. Formula (2.1) implies the first inequality below, while the second one follows from the fact that $\r^n_j$ and 
$\r^n_i$ are probability densities. 
$$d_1(\hat R^m_i(-m,s),\hat R^m_j(-m,s))\le
\sqrt{p}  ||
\frac{1}{n-1}\sum_{l\not=i}\r_l^n(-m,s,\cdot)-
\frac{1}{n-1}\sum_{l\not=j}\r_l^n(-m,s,\cdot)
||_{L^1(\T^p)}=$$
$$\frac{\sqrt{p}}{n-1}   ||
\r_j^n(-m,s,\cdot)-\r_i^n(-m,s,\cdot)
||_{L^1(\T^p)}   \le \frac{2}{n-1}  \qquad\forall s\in[-m,0],\quad
\forall i,j\in (1,\dots,n)  .  \eqno (2.14)  $$

By (2.8),
$$L^n_{c,Y^n,i}=\2|\dot x|^2-\inn{c}{\dot x}-V(t,x)-
\frac{n-1}{n}\int_{\T^p}W(x-y)\hat\r^n_i(t,y)\dr y  .  $$
By lemma 1.1, we get the second inequality below. 
$$\Big\vert\Big\vert
V(t,x)+\frac{n-1}{n}\int_{\T^p}
W(x-y)\hat\r^n_i(t,y)\dr y
\Big\vert\Big\vert_{
C((-m,0),C^3(\T^p))
}         \le$$
$$ ||
V(t,x)
||    +
||
\int_{\T^p}W(x-y)\hat\r^n_i(t,y)\dr y
||_{
C((-m,0),C^3(\T^p))
}        \le C_1  .$$
As a result, the value function $u^n_i$ satisfies point 2) of corollary 2.2; thus, $Y^n_i\in Vect_m$ and we can apply lemma 2.3, getting that $R_i^n$ belongs to $Den_m^{reg}$. Since this set is convex, (2.13) implies  that 
$\hat R_i^n\in Den^{reg}_m$; by lemma 2.4, we have that $Den_m^{reg}$ is a compact set; thus, fixing $i=1$, there is 
$n_k\tends+\infty$ such that $\hat R^{n_k}_1$ converges to 
$R\in Den^{reg}_m$; in particular, $R$ and its density $\r$ satisfy (2.5). This gives convergence only for $\hat R_1^{n_k}$; however, from (2.14) we get that 
$$\sup_{i\in(1,\dots,n_k)}\sup_{s\in[-m,0]}
d_1(\hat R^{n_k}_i(-m,s),\hat R^{n_k}_1(-m,s))\tends 0
\txt{as}k\tends+\infty    \eqno (2.15)$$
which implies that all $\hat R^{n_k}_i$ converge to the same limit $R$. 
By the same argument of (2.14), 
$$d_1(R^{n_k}(-m,s),\hat R^{n_k}_i(-m,s))\le\frac{1}{n}  .  $$
Thus, (2.15) implies the first formula of (2.11). 

\noindent{\bf Step 2.} We prove the convegence of the solutions of Hamilton-Jacobi. We set
$$\left\{
\matrix{
&W_i^{n_k}(s,x)\colon=
\frac{n-1}{n}\int_{\T^p} W(x-y)\hat\r_i^{n_k}(-m,s,y)\dr y=
\frac{n-1}{n}\int_{\T^p}W(x-y)\dr\hat R^{n_k}_i(-m,s)(y)\cr
&{}\cr
&\tilde W(s,x)\colon=\int_{\T^p} W(x-y)\r(-m,s,y)\dr y
=\int_{\T^p}W(x-y)\dr R(-m,s)(y)  . 
}    \right.
\eqno (2.16)$$
By (2.15) and the fact that $d_1$ induces weak$\ast$ convergence, we get that
$$\sup_{i\in(1,\dots,n_k)}||
W_i^{n_k}-\tilde W
||_{C([-m,0],C^3(\T^p))}\tends 0
\txt{as}k\tends+\infty . \eqno (2.17)$$
Now, $u_i^{n_k}$ is the value function of $L_{c,Y^{n_k},i}$, whose potential is $V(t,x)+W_i^{n_k}(t,x)$; by the last formula, we can apply point 3) of lemma 2.1 and get that $u_i^{n_k}$ satisfies the limit in the second formula of (2.11), with $u$ a solution of 
$(HJ)^{V+\tilde W}$ or, which is the same, of $(HJ)_{\r,f}$. 

\noindent{\bf Step 3.} We prove that the limit density $\r$ solves 
$(FP)_{-m,c-\partial_xu,\m}$.

From now on, for ease of notation, we drop the $n_k$ of the subsequence. We recall that each $\r_i^n$ solves 
$(FP)_{-m,Y^n_i,\d_{z_i}}$; by the third formula of (2.11) and lemma 2.5, we get that, if $\bar\r_i$ is a solution of 
$(FP)_{-m,c-\partial_x u,\d_{z_i}}$, then 
$$\sup_i||
\r_i^n\L^p-\bar\r_i\L^p
||_{C^0([-m,0],\Mt)}\tends 0
\txt{as}n\tends+\infty  .  $$
Now (2.2) implies the inequality below, and the last formula implies the limit. 
$$||
\r^n\L^p-\frac{1}{n}\sum_{i=1}^n\bar\r_i\L^p
||_{C([-m,0],{\cal M}_1(\T^p))}=
||
\frac{1}{n}\sum_{i=1}^n\r_i^n\L^p-
\frac{1}{n}\sum_{i=1}^n\bar\r_i\L^p
||_{C([-m,0],{\cal M}_1(\T^p))}\le$$
$$\frac{1}{n}\sum_{i=1}^n||
\r_i^n\L^p-\bar\r_i\L^p
||_{C([-m,0],{\cal M}_1(\T^p))}   \tends 0   .   $$
This means that $\r^n\L^p$ and 
$\frac{1}{n}\sum_{i=1}^n\bar\r_i\L^p$ have the same limit; we saw in step 1 that $\r^n\L^p$ converges to 
$\r\L^p$; thus, to prove that $\r\L^p$ solves 
$(FP)_{-m,c-\partial_x u, \mu}$, it suffices to prove that the limit of 
$\frac{1}{n}\sum_{i=1}^n\bar\r_i\L^p$ solves the same equation. This follows easily, since by definition 
$\frac{1}{n}\sum_{i=1}^n\bar\r_i\L^p$ solves the Fokker-Planck equation with drift $c-\partial_x u$ and initial condition 
$\frac{1}{n}(\d_{z_1}+\dots+\d_{z_n})$, and (2.10) holds.

\noindent{\bf Step 4.} We prove the last assertion of the lemma; the equality below comes from (2.6) and the fact that $Y^n$ is minimal. 
$$U^n(-m,(z_1,\dots,z_n))=
E_{w_1,\dots,w_n}\left\{
\int_{-m}^0L^n_c(s,X^n(-m,s,x),Y^n(s,X^n(-m,s,x)))\dr s
\right\}   +U(\r^n(-m,0,\cdot))  .  $$
We recall that, by corollary 2.7, 
$$Y^n=(c-\partial_x u^n_1,\dots,c-\partial_x u^n_n) . $$
Now $\r^n_i\L^p$ is the push-forward of the Wiener measure by $X^n_i$, and the Brownian motions $w_i$ are independent. This implies that
$$U^n(-m,(z_1,\dots,z_n))=
\frac{1}{n}\sum_{i=1}^n
\int_{[-m,0]\times\T^p}L_c(t,x,c-\partial_xu^n_i)\r^n_i(t,x)\dt\dx-  $$
$$\frac{1}{2n^2}\sum_{i\not=j\in(1,\dots,n)}\int_{[-m,0]\times\T^p\times\T^p}
W(x_i-x_j)\r_i^n(-m,t,x_i)\r_j^n(-m,t,x_j)\dr x_i\dr x_j\dt
+\int_{\T^p}f(x)\r^n(-m,0,x)\dx   .   $$
Using (2.11), we get immediately that
$$U^n(-m,(z_1,\dots,z_n))\tends
\int_{[-m,0]\times\T^p}\L_{c,\2\r}(t,x,c-\partial_xu)\r(-m,0,x)\dt\dx
+\int_{\T^p}f(x)\r(-m,0,x)\dx  .  $$

\fin

\noindent{\bf Proof of theorem 3.} Let the measure $\m$, the couple $(u,\r)$ and the function $U(-m,\mu)$ be as in the last lemma; let the operator $\Lambda^m_c$ be as in the introduction, and let $Y=c-\partial_xu$. We are going to prove that 
$$(\Lambda^m_cU)(\mu)=
U(-m,\m)=E_w\left\{
\int_{-m}^0\L_{c,\2\r}(s, X, Y)\dr s+f(X(-m,0,\mu))
\right\}    =$$
$$\min_{\tilde Y} E_w\left\{
\int_{-m}^0\L_{c,\2\tilde\r}(s,\tilde X,\tilde Y)\dr s
+f(\tilde X(-m,0,\mu))
\right\}  .  \eqno (2.18)$$
The functions $\r$ and $\tilde\r$ in the formula above satisfy 
$(FP)_{-m,Y,\m}$ and $(FP)_{-m,\tilde Y,\m}$ respectively, while $X$ and $\tilde X$ satisfy $(SDE)_{-m,Y,\mu}$ and 
$(SDE)_{-m,\tilde Y,\mu}$ respectively. The minimum is taken over all Lipschitz vector fields $\tilde Y$.

Note that, in principle, $U(-m,\mu)$ could depend on the subsequence $\{ n_k \}_{k\ge 1}$ chosen in lemma 2.8; the formula above says that this is not the case. Moreover, it says that any $(u,\r)$ arising in lemma 2.8 as the limit of a subsequence, minimizes the last expression of (2.18).

The second equality of (2.18) follows from lemma 2.8: it is just another way of writing (2.12). Again by lemma 2.8, 
$(u,\r)\in{\bf S}$, and thus, by the definition of 
$(\Lambda^m_c U)(\mu)$, 
$$(\Lambda^m_c U)(\mu)\le 
E_w\left\{
\int_{-m}^0\L_{c,\2\r}(s, X, Y)\dr s
+f( X(-m,0,\mu))
\right\}=
U(-m,\mu)  .  \eqno (2.19)$$
Now we prove that 
$$U(-m,\mu)\le
\inf_{\tilde Y}
E_w\left\{
\int_{-m}^0\L_{c,\2\tilde\r}(s,\tilde X,\tilde Y)\dr s
+f(\tilde X(-m,0,z))
\right\}   .  \eqno (2.20)$$
To prove this, we consider the n-particle value function 
$U^n(-m,z_1,\dots,z_n)$. Let $\tilde Y$ be a Lipschitz vector field on $[-m,0]\times\T^p$. Let $\tilde X^n_i$ solve $(SDE)_{-m,\tilde Y,\d_{z_i}}$; let us suppose that (2.10) holds. Let us set 
$\tilde X^n=(\tilde X_1^n,\dots,\tilde X_n^n)$ and 
$\tilde Y^n=(\tilde Y,\dots,\tilde Y)$. Let 
$\tilde\r$ solve $(FP)_{-m,\tilde Y,\mu}$ and  let $\tilde\r_i$ solve 
$(FP)_{-m,\tilde Y,\delta_{z_i}}$; by linearity, we get that 
$\frac{1}{n}\sum_{i=1}^n\tilde\r_i(-m,s,x)$ solves 
$(FP)_{-m,\tilde Y,\frac{1}{n}(\d_{z_1+\dots+\d_{z_n}})}$. In other words, $\tilde\r$ and $\frac{1}{n}\sum_{i=1}^n\tilde\r_i(-m,s,x)$ solve a Fokker-Planck equation with the same drift, but initial distributions $\mu$ and $\frac{1}{n}(\d_{z_1}+\dots\d_{z_n})$ respectively; by (2.10), it is standard to see that
$$\sup_{s\in[-m,0]}   d_1\left(
\frac{1}{n}\sum_{i=1}^n\tilde\r_i(-m,s,x)\L^p,\tilde\r(-m,s,x)\L^p
\right)
\tends 0   \txt{as}n\tends+\infty  .   \eqno (2.21)$$
We set $z=(z_1,\dots,z_n)$ and
by (2.6) we get the inequality below; the first equality is the definition of $L^n_c$, the second one comes from the fact that the Brownian motions 
$w_1,\dots,w_n$ are independent.
$$U^n(-m,z)\le$$ 
$$E_{w_1,\dots,w_n}\left\{
\int_{-m}^0
L^n_c(s,\tilde X^n(-m,s,z),\tilde Y^n(s,\tilde X^n(-m,s,z)))\dr s
+\frac{1}{n}\sum_{i=1}^nf(\tilde X_i(-m,0,z_i))
\right\}  =$$
$$E_{w_1,\dots,w_n}\left\{
\frac{1}{n}\sum_{i=1}^n
\int_{-m}^0L_c(s,\tilde X_i,\tilde Y)\dr s-
\frac{1}{2n^2}\sum_{i\not=j}
\int_{-m}^0W(\tilde X_i-\tilde X_j)\dr s+
\frac{1}{n}\sum_{i=1}^n f(\tilde X_i(-m,0,z_i))
\right\}    =   $$
$$\frac{1}{n}\sum_{i=1}^n
\int_{[-m,0]\times\T^p}L_c(s,x,\tilde Y(s,x))\tilde\r_i(-m,s,x)\dr s\dx-$$
$$\frac{1}{2n^2}\sum_{i\not=j=1}^n
\int_{[-m,0]\times\T^p\times\T^p}
W(x_i-x_j)\tilde\r_i(-m,s,x_i)\tilde\r_j(-m,s,x_j)\dr s\dx_i\dx_j+
\frac{1}{n}\sum_{i=1}^n 
\int_{\T^p}f(x)\tilde\r_i(-m,0,x)\dx  .  $$
We take limits in the formula above, using the last assertion of lemma 2.8 for the left hand side and (2.21) for the right hand side; we get that
$$U(-m,\mu)\le
\int_{-m}^0\L_{c,\2\tilde\r}(s,x,\tilde Y)\tilde\r(t,x)\dt\dx 
+\int_{\T^p}f(x)\tilde\r(-m,0,x)\dx .  $$
Since $\tilde Y$ is an arbitrary Lipschitz vector field, we get that
(2.20) holds.

Let now $(\bar u,\bar\r)\in{\bf S}$ be minimal in the definition of 
$(\Lambda^m_c U)(\mu)$; setting 
$\tilde Y=c-\partial_x\bar u$, (2.20) implies the inequality below, while the equality comes from our choice of $\tilde Y$.
$$U(-m,\mu)\le
E_w\left\{
\int_{-m}^0\L_{c,\2\tilde\r}(s,\tilde X,\tilde Y)\dr s
+f(\tilde X(-m,0,z))
\right\}    =(\Lambda^m_c U)(\mu)  .  $$
This yields the inequality opposite to (2.19). In other words, we have proven the first equality of (2.18); the second one, as we have seen, is lemma 2.8. As for the third one, it suffices to prove the inequality opposite to (2.20), which we do presently.

Let $(z_1,\dots,z_n)$ satisfy (2.10), let 
$Y^n=(Y_1^n,\dots,Y^n_n)$ be minimal in (2.6), and let us set 
$\tilde Y^n=Y^n_1$. Let $\tilde\r^n$ satisfy $(FP)_{-m,\tilde Y^n,\m}$. By (2.11) and (2.17), we obtain that there is $\g_n\tends 0$ such that
$$E_w\left\{
\int_{-m}^0\L_{c,\2\tilde\r^n}(s,\tilde X^n,\tilde Y^n)\dr s
+f(\tilde X^n(-m,0,\mu))
\right\}    \le$$
$$E_{w_1,\dots,w_n}\left\{
\int_{-m}^0[
\frac{1}{n}\sum_{i=1}^nL_c(t,X^n_i,Y^n_i)+
\frac{1}{2n^2}\sum_{i,j=1}^n W(X^n_i-X^n_j)
]\dr s    
+\frac{1}{n}\sum_{i=1}^n f(X^n_i(-m,0,\mu))
\right\}
+\g_n   .   $$
Since the limit of the term on the right is $U(-m,\mu)$ by lemma 2.8, we get that
$$\inf_{\tilde Y}\left\{
\int_{-m}^0\L_{c,\2\tilde\r}(s,\tilde X,\tilde Y)\dr s
+f(\tilde X^n(-m,0,\mu))
\right\}    \le U(-m,\m)  $$
yielding the inequality opposite to (2.20).

\fin

We need the following lemma to prove the semigroup property.

\lem{2.9} Let $Y_1$ be a Lipschitz vector field on 
$[-(n+m),-n]\times\T^p$, and let $Y_2$ be a Lipschitz vector field on $[-n,0]\times\T^p$. Then, for all $\e,\d\in(0,1)$, there is a Lipschitz vector field $Y$ which coincides with $Y_1$ when 
$t\in[-(n+m),-n]$, and with $Y_2$ when $t\in[-n+\d,0]$. Moreover, 
$Y$ satisfies 
$$E_w\left\{
\int_{-(n+m)}^{0}
\L_{c,\2\r_Y}(s,X(-(n+m),s,\m),Y(s,X(-(n+m),s,\m)))\dr s
\right\}  \le$$
$$E_w\left\{
\int_{-(n+m)}^{-n}
\L_{c,\2\r_{Y_1}}(s,X_1(-(n+m),s,\m),
Y_1(s,X_1(-(n+m),s,\m)))\dr s
\right\}  +$$
$$E_w\left\{
\int_{-n}^{0}
\L_{c,\2\r_{Y_2}}(s,X_2(-n,s,\r_{Y_1}(-n)),
Y_2(s,X_2(-n,s,\r_{Y_1}(-n))))\dr s 
\right\} +\e  .  \eqno (2.22)$$
In the formula above, $X_1$ solves $(SDE)_{-(n+m),Y_1,\mu}$, 
$X_2$ solves $(SDE)_{-m,Y_2,\r_1(-m)}$, $X$ solves
 
\noindent $(SDE)_{-(n+m),Y,\mu}$ and $\r_1$, $\r_2$, $\r_Y$ are the densities of the laws of $X_1$, $X_2$ and $X$ respectively. 

Moreover, we can require that 
$$|
U(\r_Y(0)\L^p)-U(\r_{Y_2}(0)\L^p)
|    \le\e   .  \eqno (2.23)   $$

\proof Let $\bar\d\in(0,\d)$; it is always possible to find a Lipschitz vector field $Y$ coinciding with $Y_1$ on $[-(n+m),-n]\times\T^p$ and with $Y_2$ on $[-n+\bar\d,0]\times\T^p$, and such that 
$||Y||_\infty$ is bounded uniformly in $\bar\d$; we forego the easy proof of this fact. 

We note that $X=X_1$ when 
$s\in[-(n+m),-n]$, since both functions solve the same stochastic differential equation; as a consequence,
$$E_w\left\{
\int_{-(n+m)}^{-n}
\L_{c,\2\r_Y}(s,X(-(n+m),s,\m),Y(s,X(-(n+m),s,\m)))\dr s
\right\}  =$$
$$E_w\left\{
\int_{-(n+m)}^{-n}
\L_{c,\2\r_{Y_1}}(s,X_1(-(n+m),s,\m),
Y_1(s,X_1(-(n+m),s,\m)))\dr s
\right\}  .  $$
Thus, it suffices to prove that
$$E_w\left\{
\int_{-n}^{0}
\L_{c,\2\r_Y}(s,X(-n,s,\r_{Y_1}(-n)),Y(s,X(-n,s,\r_{Y_1}(-n))))
\dr s
\right\}  \le$$
$$E_w\left\{
\int_{-n}^{0}
\L_{c,\2\r_{Y_2}}(s,X_2(-n,s,\r_{Y_1}(-n)),
Y_2(s,X_2(-n,s,\r_{Y_1}(-n))))\dr s 
\right\} +\e  .  \eqno (2.24)$$
To prove this, we recall that $X$ and $X_2$ solve two stochastic differential equations with drift $Y$ and $Y_2$ respectively; this means that, for $s\ge -n$ and any trajectory $w$ of the Brownian motion, we have that
$$X(-n,s,\r_{Y_1}(-n))(w)=X(-n,-n,\r_{Y_1}(-n))(w)+
\int_{-n}^sY(\tau,X(-n,\tau,\r_{Y_1}(-n))(w))\dr\tau+
w(s)-w(-n)     \eqno (2.25)     $$
and 
$$X_2(-n,s,\r_{Y_1}(-n))(w)=X_2(-n,-n,\r_{Y_1}(-n))(w)+
\int_{-n}^sY_2(\tau,X_2(-n,\tau,\r_{Y_1}(-n))(w))\dr\tau+
w(s)-w(-n)  .  \eqno (2.26)    $$
Since $X_2(-n,-n,\r_{Y_1}(-n))$ and $X(-n,-n,\r_{Y_1}(-n))$ have the same law $\r_{Y_1}(-n)$, we can as well suppose that
$$X_2(-n,-n,\r_{Y_1}(-n))(w)=
X(-n,-n,\r_{Y_1}(-n))(w)      \eqno (2.27)$$
for all realizations $w$ of the Brownian motion. 
Subtracting (2.26) from (2.25) and using the formula above, we get the inequality below; the equality is the definition of the function $a$. 
$$|
X(-n,s,\r_{Y_1}(-n))(w)-X_2(-n,s,\r_{Y_1}(-n))(w)
|\le$$
$$\int_{-n}^s|
Y(\tau,X(-n,\tau,\r_{Y_1}(-n))(w))-
Y_2(\tau,X_2(-n,\tau,\r_{Y_1}(-n))(w))
|   \dr\tau=$$
$$\int_{-n}^s a[\tau,X(-n,\tau,\r_{Y_1}(-n))(w),
X_2(-n,\tau,\r_{Y_1}(-n))(w)]\dr\tau   $$
Since $Y$ and $Y_2$ are  bounded uniformly in $\bar\d$, we get that $|a|\le M$ if $\tau\in[-n,0]$ for a constant $M$ independent on $\bar\d$; since $Y$ coincides with the Lipschitz 
$Y_2$ on $[-n+\bar\d,0]$, we get that, for $\tau\ge -n+\bar\d$,
$$|
a(\tau,x,y)
|\le K|x-y|$$
for a constant $K$ independent on $\bar\d$. From the last two formulas, we get that
$$|
X(-n,s,\r_{Y_1}(-n))(w)-X_2(-n,s,\r_{Y_1}(-n))(w)
|\le$$
$$\int_{-n}^{-n+\bar\d}M\dr\tau+
\int_{-n+\bar\d}^sK|
X(-n,\tau,\r_{Y_1}(-n))(w)-X_2(-n,\tau,\r_{Y_1}(-n))(w)
|    \dr\tau     .    $$
Using the Gronwall lemma and (2.27), we get that there is a function 
$\g(\bar\d)$, tending to zero as $\bar\d$ tends to zero, such that
$$|
X(-n,s,\r_{Y_1}(-n))(w)-\tilde X(-n,s,\r_{Y_1}(-n))(w)
|   \le\g(\bar\d)  $$
for all realizations $w$ of the Brownian motion. From this, (2.24) follows easily.

On the other hand, it is easy to see that the formula above implies that, as $\bar\d\tends 0$, $\r_Y(0)\L^p$ converges weak$\ast$ to 
$\r_{Y_2}(0)\L^p$. Since $U$ is Lipschitz for the 1-Wasserstein distance, (2.23) follows.

\fin

\prop{2.10} 1) The map $\Psi^m_c$ defined in the introduction has the semigroup property, i. e. for $n,m\ge 0$ and 
$U\in C({\cal M}_1(\T^p),\R)$,
$$\Psi_c^{n+m}U=\Psi_c^n\circ\Psi_c^mU  .  $$

\noindent 2) If $U\le V\in C({\cal M}_1(\T^p),\R)$, then 
$\Psi^m_cU\le\Psi^m_cV$.

\noindent 3) For all $a\in\R$ and $U\in C({\cal M}_1(\T^p),\R)$, 
$\Psi^m_c(U+a)=(\Psi^m_cU)+a$.

\proof Properties 2) and 3) follow in a standard way from the definition of $\Psi^m_c$; we prove 1).

Let $\m\in{\cal M}_1$; by the definition of 
$\Psi_c^{n+m}U$ as an infimum, for any $\e>0$ we can find a Lipschitz vector field $Y$ such that
$$\Psi_c^{n+m}U(\m)\ge
E_w\left\{
\int_{-(n+m)}^0\L_{c,\2\r_Y}(s,X(-(n+m),s,\m),Y(s,X(-(n+m),s,\m)))
\dr s
\right\}      + U(\r_Y(0)\L^p)  -\e   $$
where $X$ solves $(SDE)_{-(n+m),Y,\m}$ and $\r_Y$ is, as usual, the solution of the Fokker-Planck equation with initial condition 
$\m$. By the Chapman-Kolmogorov formula, the formula above implies the first inequality below. 
$$(\Psi^{n+m}_cU)(\m)\ge E_w\left\{
\int_{-(n+m)}^{-n}
\L_{c,\2\r_Y}(s,X(-(n+m),s,\m),Y(s,X(-(n+m),s,\m)))\dr s
\right\}  +$$
$$E_w\left\{
\int_{-n}^0\L_{c,\2\r_Y}(s,X(-n,s,\r_Y(-n)),Y(s,X(-n,s,\r_Y(-n))))\dr s
\right\}  +U(\r_Y(0)\L^p)-\e\ge$$
$$E_w\left\{
\int_{-(n+m)}^{-n}
\L_{c,\2\r_Y}(s,X(-(n+m),s,\m),Y(s,X(-(n+m),s,\m)))\dr s
\right\}  +    (\Psi^n_cU)(\r_Y(-n)\L^p)-\e\ge$$
$$(\Psi^m_c\circ\Psi^n_c U)(\mu)-\e  .  $$
The second and third inequalities above come from the definition of $\Psi^n_cU$ and $\Psi^m_c\circ\Psi^n_cU$ as infima. Since 
$\e$ is arbitrary, this means that 
$$(\Psi^{n+m}_cU)(\mu)\ge
(\Psi^n_c\circ\Psi^n_c U)(\mu)  .  \eqno (2.28)$$

We prove the opposite inequality. By the definition of 
$\Psi^m_c\circ\Psi^n_c (U)$, we can find a Lipschitz vector field $Y_1$ such that
$$(\Psi_c^{m}\circ\Psi_c^n U)(\m)\ge
E_w\left\{
\int_{-(m+n)}^{-n}
\L_{c,\2\r_{Y_1}}(s,X_1(-(n+m),s,\m),Y_1(s,X_1(-(n+m),s,\m)))
\dr s
\right\}      +$$
$$\Psi_c^n( U) (\r_{Y_1}(-n)\L^p)  -\e  .  \eqno (2.29)$$
By the definition of $\Psi^n_cU$, we can find another Lipschitz vector field $Y_2$ such that
$$\Psi_c^{n}U(\r_{Y_1}(-n))\ge
E_w\left\{
\int_{-n}^{0}
\L_{c,\2\r_{Y_2}}(s,X_2(-(n+m),s,\r_{Y_1}(-n)),
Y_2(s,X_2(-(n+m),s,\r_{Y_1}(-n))))
\dr s
\right\}      +$$
$$ U(\r_{Y_2}(0)\L^p)  -\e  .  \eqno (2.30)$$
Let $\e,\d>0$; by lemma 2.9, we can find a Lipschitz vector field $Y$ equal to $Y_1$ on $[-(n+m),-n]\times\T^p$ and to $Y_2$ on 
$[-n+\d,0]\times\T^p$, such that (2.22) and (2.23) holds.  The first inequality below comes from the definition of $\Psi^{n+m}U$ as a infimum; the second one from (2.22) and (2.23); the third and fourth ones come from (2.30) and (2.29) respectively.
$$(\Psi_c^{n+m}U)(\m)\le
E_w\left\{
\int_{-(n+m)}^0\L_{c,\2\r_Y}(s,X(-(n+m),s,\m),Y(s,X(-(n+m),s,\m)))\dr s
\right\}  +  U(\r_Y(0)\L^p)\le$$
$$E_w\left\{
\int_{-(n+m)}^{-n}
\L_{c,\2\r_{Y_1}}(s,X_1(-(n+m),s,\m),Y_1(s,X_1(-(n+m),s,\m)))\dr s
\right\}   +$$
$$E_w\left\{
\int^{0}_{-n}
\L_{c,\2\r_{Y_2}}(s,X_2(-n,s,\r_{Y_1}(-n)),Y_2(s,X_2(-n,s,\r_{Y_1}(-n))))\dr s
\right\}   + U(\r_{Y_2}(0)\L^p)+2\e\le$$
$$E_w\left\{
\int_{-(n+m)}^{-n}
\L_{c,\2\r_{Y_1}}(s,X_1(-(n+m),s,\m),Y_1(s,X_1(-(n+m),s,\m)))\dr s
\right\}   + (\Psi^n_cU)(\r_{Y_1}(-n)\L^p)+3\e\le$$
$$(\Psi^n_c\circ\Psi^m_cU)(\m)  +4\e .  $$
Since $\e>0$ is arbitrary, we get the inequality opposite to (2.28), and thus the thesis.

\fin

\vskip 2pc
\centerline{\bf\S 3}
\centerline{\bf Fixed points}
\vskip 1pc

As in [8] and in [15], the following proposition is essential in proving theorem 4.

\prop{3.1} Let $U$ be linear as in theorem 2. Then, there is $L>0$, independent on $n$, such that $\Lambda^n_c U=\Psi^n_c U$ is $L$-Lipschitz for the Wasserstein distance $d_1$.

\rm
\vskip 1pc

To prove this proposition, we shall need two lemmas. 

\lem{3.2} Let $R\in C([-m,0],\Mt)$ and let $u$ solve $(HJ)_{R,f}$. Then, there is $C>0$, independent both of $m\in\N$ and of 
$R\in C([-m,0],\Mt)$, such that
$$||\partial_xu(t,\cdot)||_{C^1([-m,0],C(\T^p))}+
||\partial_xu(t,\cdot)||_{C([-m,0],C^2(\T^p))}\le C  .  $$

\proof We have seen in section 1 that, if $u$ is a solution of 
$(HJ)_{R,f}$, $v=e^{-\b u}$ and $a\in\R$, then 
$e^{-\b a}v=e^{-\b(u+a)}$ is a solution of 
$(TS)_{R,e^{-\b(f+a)}}$ with $A=0$. Let $a_k\in\R$ such that 
$e^{-\b a_k}v(-k,\cdot)=e^{-\b(u+a_k)}$ satisfies (1.7) for 
$k=0,1,2,\dots$. By the Feynman-Kac formula, for $k\ge 0$,
$$e^{-\b a_k}v(-k-1,\cdot)=
L_{(\psi,0,-1)}(
e^{-\b a_k}v(-k,\cdot)
)   .  \eqno (3.1)$$
Since $e^{-\b a_k}v(-k,\cdot)$ satisfies (1.7), formulas (1.10) and (1.11) hold and we get that, for $k\ge 0$,
$$\frac{1}{C_1}\le e^{-\b a_k}v(-k-1,x)\le C_1  
	\quad\forall x\in\T^p  .  $$
We consider (3.1) with $e^{-\b a_k}v(-k-1,x)$ on the right hand side and differentiate under the integral sign; proceeding as in lemma 1.3, and using the last formula, we get that, for $k\ge 0$, 
$$||
e^{-\b a_k}v(-k-2,\cdot)||_{C^3(\T^p)}\le C_2     .    $$
As in lemma 1.4, this implies that there is $C_3>0$, independent on $k\ge 0$ (it depends only on $C_1$ and $C_2$) such that, for 
$k\ge 0$,
$$\txt{if} t\in[-(k+3),-(k+2)], \txt{then}
\frac{1}{C_3}\le e^{-\b a_{k}}v(t,\cdot)\le C_3 \txt{and}  $$
$$||e^{-\b a_{k}}v(t,\cdot)||_{C([-(k+3),-(k+2)],C^3(\T^p))}+
||e^{-\b a_{k}}v(t,\cdot)||_{C^1([-(k+3),-(k+2)],C^1(\T^p))}\le C_3    .  $$
By our definition of $v$,
$$\txt{for}t\in[-(k+3),-(k+2)],\qquad
u=-\frac{1}{\b}\log(e^{-\b a_{k}}v)-a_{k}  .  $$
From the two formulas above, we get that
$$||
\partial_xu
||_{C([-m,-2],C^2(\T^p))}+
||
\partial_xu
||_{C^1([-m,-2],C(\T^p))}\le C    \txt{for}  t\le -2 . $$
It remains to bound $\partial_xu(t,x)$ when $t\in[-2,0]$; since 
$f\in C^3(\T^p)$, this follows by differentiation under the integral sign in (1.20), and we are done.

\fin

We recall some notation: in the following $U^n(-m,z)$ will be the minimum in (2.6); moreover, given 
$z=(z_1,\dots,z_n)\in(\R^{p})^n$, we set
$$|z|_1=|z_1|+\dots+|z_n|    $$
and we define $z^\prime\in(\R^p)^{n-1}$ by 
$z=(z_1,\dots,z_n)=(z_1,z^\prime)$.

\lem{3.3} Let $U$ be as in theorem 2 and let $U^n(-m,z)$ be defined as in (2.6). Then, there is a constant $C>0$ such that, for all positive integers $n$ and $m$, we have
$$|
U^n(-m,z)-U^n(-m,\tilde z)
|   \le
\frac{C}{n}|z-\tilde z|_1 . $$

\proof It suffices to prove that, for $i=1,\dots,n$, the function 
$$\fun{}{z_i}{U^n(-m,z_1,\dots,z_i,\dots,z_n)}$$
is $\frac{C}{n}$-Lipschitz and that the constant $C$ does not depend neither on $m\in\N$ nor on 
$(z_1,\dots,z_{i-1},z_{i+1},\dots,z_n)\in(\T^p)^{n-1}$. We shall prove this for $i=1$; from this the general case follows, since 
$U$ is a symmetric function of $(z_1,\dots,z_n)$.

We write the function $U^n(-m,(z_1,z^\prime))$ as in lemma 2.6,  isolating particle $z_1$:
$$U^n(-m,(z_1,z^\prime))=
E_{w_1,\dots,w_n}\bigg\{
\frac{1}{n}\sum_{j\not =1}
\int_{-m}^0L_c(s,X_j(-m,s,z_j),Y_j(s,X_j(-m,s,z_j)))\dr s+$$
$$\frac{1}{2n^2}\sum_{j,i\not=1}
\int^0_{-m}W(X_i(-m,s,z_i)-X_j(-m,s,z_j))\dr s
+\frac{1}{n}\sum_{j\not=1}f(X_j(-m,0,z_j))
\bigg\}  +$$
$$\frac{1}{n}E_{w_1}\left\{
\int_{-m}^0L^n_{c,Y,1}(s,X_1(-m,s,z_1),
Y_1(s,X_1(-m,s,z_1)))\dr s
+f(X_1(-m,0,z_1))
\right\}   \eqno (3.2)  $$
where the vector field $Y=(Y_1,\dots,Y_n)$ is minimal. By the definition of $U^n(-m,(\tilde z_1,z^\prime))$ as a minimum, we get that
$$U^n(-m,(\tilde z_1,z^\prime))\le
E_{w_1,\dots,w_n}\bigg\{
\frac{1}{n}\sum_{j\not =1}
\int_{-m}^0L_c(s,X_j(-m,s,z_j),Y_j(s,X_j(-m,s,z_j)))\dr s+$$
$$\frac{1}{2n^2}\sum_{j,i\not=1}
\int^0_{-m}W(X_i(-m,s,z_i)-X_j(-m,s,z_j))\dr s
+\frac{1}{n}\sum_{j\not=1}f(X_j(-m,0,z_j))
\bigg\}  +$$
$$\frac{1}{n}E_{w_1}\left\{
\int_{-m}^0L^n_{c,Y,1}(s,X_1(-m,s,\tilde z_1),
Y_1(s,X_1(-m,s,\tilde z_1)))\dr s
+f(X_1(-m,0,\tilde z_1))
\right\}   .   \eqno (3.3)$$
The term with $E_{w_1,\dots,w_n}$ is identical in (3.2) and (3.3); defining $u^n_1$ as in corollary 2.7, we see that
the term with $E_{w_1}$ is equal to the function $u^n_1(-m,z_1)$ in (3.2), and to $u^n_1(-m,\tilde z_1)$ in (3.3); subtracting (3.2) from (3.3), we get that
$$U^n(-m,(\tilde z_1,z^\prime))-U^n(-m,(z_1,z^\prime))\le
\frac{1}{n}[
u^n_1(-m,\tilde z_1)-u^n_1(-m,z_1)
] \le
\frac{C}{n}|\tilde z_1-z_1|   $$
where the last inequality comes from lemma 3.2. Exchanging the 
r\^oles of $z_1$ and $\tilde z_1$, we get that the function 
$\fun{}{z_1}{U^n(-m,(z_1,z^\prime))}$ is $\frac{C}{n}$-Lipschitz; we saw at the beginning of the proof that this implies the thesis.

\fin

\noindent{\bf Proof of proposition 3.1.} By lemma 2.8, we know that
$$\txt{if}\frac{1}{n}(
\d_{z_1}+\dots+\d_{z_n}
)      \tends  \mu       \txt{then}
U^n(-m,z_1,\dots,z_n)\tends U(-m,\mu)   .   $$
We saw in (2.18) that $U(-m,\mu)=(\Lambda^m_cU)(\mu)$. Thus it  suffices to show that 
$$|
U^n(-m,(x_1,\dots,x_n))-U^n(-m,(y_1,\dots,y_n))
|   \le
Ld_1\left(
\frac{1}{n}(\d_{x_1}+\dots+\d_{x_n}),
\frac{1}{n}(\d_{y_1}+\dots+\d_{y_n})
\right)  .  $$
It is standard ([5]) that
$$d_1(
\frac{1}{n}(\d_{x_1}+\dots+\d_{x_n}),
\frac{1}{n}(\d_{y_1}+\dots+\d_{y_n})
)   =\min_\s
\frac{1}{n}\sum_{i=1}^n|x_i-y_{\s(i)}|    \eqno (3.7)$$
where the minimum is taken over all the permutations $\s$ of 
$\{ 1,\dots,n \}$. In terms of transport, when we are connecting two $n$-uples of deltas, there is not just a minimal transfer plan, but a minimal transfer map.

Since
$$U^n(-m,(y_1,\dots,y_n))=U^n(-m,(y_{\s(1)},\dots,y_{\s(n)}))  ,  $$
we have to prove that, for $\s$ minimal in (3.7), 
$$|
U^n(-m,(x_1,\dots,x_n))-U^n(-m,(y_{\s(1)},\dots,y_{\s(n)}))
|    \le \frac{C}{n}\sum_{i=1}^n|x_i-y_{\s(i)}|   .   $$
But this is an immediate consequence of lemma 3.3. 

\fin

When $U$ is linear, we define $\Lambda_{c,\l}U=\Lambda_cU+\l$; thus, in case of a linear $U$, we have that 
$\Lambda_{c,\l}U=\Psi_{c,\l}U$ for the operator $\Psi_{c,\l}$ defined in the introduction. In the next lemma, we stick to the 
$\Lambda_{c,\l}$ notation.

\lem{3.4} Let the operator $\Lambda_{c,\l}$ be defined as in the introduction and let $U=0$. Then, there is a unique $\l\in\R$ such that
$$\hat U(\mu)\colon=
\liminf_{n\tends+\infty}(\Lambda^n_{c,\l}0)(\m)  $$
is finite for all $\m\in{\cal M}_1$. Moreover, $\hat U$ is 
$L$-Lipschitz for the constant $L$ of proposition 3.1.

\proof Clearly, there is at most one $\l\in\R$ for which the $\liminf$ above is finite; let us prove that it exists. This means finding 
$\l\in\R$ such that, for all $\mu\in{\cal M}_1(\T^p)$,
$$-\infty<\liminf_{n\tends+\infty}(\Lambda^n_{c,\l}0)(\mu)<+\infty . \eqno (3.8)$$
Note that the formula above implies that $\hat U$ is finite; it is 
$L$-Lipschitz because it is the $\liminf$ of $L$-Lipschitz functions. 

By proposition 3.1, $\Lambda_{c,0}^n0$ is $L$-Lipschitz for all 
$n\in\N$; since ${\cal M}_1(\T^p)$ is a compact metric space, we can find $M>0$ such that
$$\max\Lambda_{c,0}^n0-
\min\Lambda_{c,0}^n0\le M\quad\forall n\ge 1  .  
\eqno (3.9)$$
Possibly taking a larger $M$, we can suppose that 
$$||\Lambda^1_{c,0}0||_{\sup}\le M  .  $$
By point 2) of proposition 2.10, this implies the first inequality below; the equality follows by point 1), and the second inequality by point 3) of the same proposition. 
$$(\Lambda_{c,0}^20)(\mu)=
(\Lambda^1_{c,0} (\Lambda^1_{c,0} 0))(\mu)\le 
\Lambda^1_{c,0}(0+M)\le 2M   .  $$
Exchanging signs, this implies that
$$||\Lambda_{c,0}^2||_{\sup}\le 2M   .  $$
Iterating, we get
$$||\Lambda^n_{c,0}0||_{\sup}\le n M\qquad
\forall n\ge 1  .  \eqno (3.10)$$
We set 
$$a_n=\min_\mu(\Lambda^n_{c,0}0)(\mu)\txt{and}
-\l=\liminf_{n\tends+\infty}\frac{a_n}{n}  .  \eqno (3.11)$$
From (3.10), it follows that 
$\l\in[-M,M]$. We assert that $\l$ satisfies (3.8). We prove the inequality on the left of (3.8), since the one on the right is analogous; actually, we are going to prove that, for all 
$\mu\in{\cal M}_1(\T^p)$ and all $n\in\N$,
$$(\Lambda_{c,\l}^n0)(\mu)>-10M  .  $$
Indeed, let us suppose by contradiction that, for some $m\in\N$ and $\bar\mu\in{\cal M}_1(\T^p)$, we have
$$(\Lambda^m_{c,\l}0)(\bar\mu)\le -10M  .  $$
By (3.9), this implies that, for all $\mu\in{\cal M}_1(\T^p)$,
$$(\Lambda^m_{c,\l}0)(\mu)\le -9M  .  \eqno (3.12)$$
Let $\mu\in{\cal M}_1(\T^p)$; the first inequality below comes from (3.12) and points 1) and 2) of proposition 2.10, the equality from point 3) of the same proposition, the last inequality from (3.12).
$$(\Lambda^{2m}_{c,\l}0)(\mu)\le
[
\Lambda^m_{c,\l}(-9M)
]  (\mu)=
-9M+(\Lambda^m_{c,\l}0)(\mu)\le -18M  .  $$
Proceeding by induction, we find that
$$(\Lambda^{km}_{c,\l}0)(\mu)\le  -9kM  \qquad
\forall\mu\in \Mt  .  $$
Since 
$$(\Lambda^{km}_{c,\l}0)(\mu)=(\Lambda^{km}_{c,0}0)(\mu)+
km\l  ,  $$
we get that
$$(\Lambda^{km}_{c,0}0)(\mu)\le -9kM-(km)\l   \qquad
\forall \mu\in\Mt  .  $$
By the definition of $\l$ in (3.11), this implies that 
$-\l\le-\l-\frac{9M}{m}$; this contradiction proves (3.8) and thus the lemma.

\fin

\noindent{\bf Proof of theorem 4.} Let $\hat U$ be as in the last lemma; since $\hat U$ may not be linear, we switch to the $\Psi^1_{c,\l}$ notation. We have to prove that $\hat U$ is a fixed point of 
$\Psi^1_{c,\l}$ and that a minimizing vector field exists.

Let $\mu\in{\cal M}_1(\T^p)$ and $\e>0$; by the definition of $\Psi^1_{c,\l}\hat U$, we can find a Lipschitz vector field 
$\bar Y$ for which
$$E_w\left\{
\int_{0}^{1}\L_{c,\2\bar\r}(s,\bar X,\bar Y)\dr s
\right\}   
+\hat U(\bar\r(1)\L^p)+\l\le (\Psi^1_{c,\l}\hat U)(\mu)+\e  .   
\eqno (3.13) $$
To use this formula, we are going to express 
$\hat U(\bar\r(1)\L^p)$ by the limit of lemma 3.4.
 
Let  $n\in\N$; by theorem 3, applied to $f\equiv 0$ with an obvious translation in time, we can find $Y_n$ be such that
$$(\Psi^n_{c,0}0)(\mu)=
\min_Y E_w\left\{
\int_{1}^{n+1}\L_{c,\2\r}(s,X,Y)
\right\}  =
E_w\left\{
\int_{1}^{n+1}\L_{c,\2\r_n}(s,X_n,Y_n)
\right\}     $$
where $\r_n$ stays for $\r_{Y_n}$ and $\r$ for $\r_Y$; the initial time for $(SDE)$ and $(FP)$ is $1$. 
Let $\tilde Y$ be equal to $\bar Y$ on $[0,1]\times\T^p$, to $Y_n$ on $[1+\d,n+1]\times\T^p$ and a Lipschitz connection in between. By lemma 2.9, we can choose the Lipschitz connection in such a way that
$$E_w\left\{
\int_{0}^{1} 
\L_{c,\2\bar\r}(s,\bar X(0,s,\mu),\bar Y(s,\bar X(0,s,\mu)))
\dr s
\right\}  +  $$
$$E_w\left\{
\int_{1}^{n+1} 
\L_{c,\2\r_n}(s, X_n(1,s,\mu),Y_n(s,X_n(1,s,\r_{\bar Y}(1))))
\dr s
\right\}  
\ge$$
$$
E_w\left\{
\int_{0}^{n+1}
\L_{c,\2\tilde\r}(s,\tilde X(0,s,\mu),\tilde Y(s,\tilde X(0,s,\mu)))
\dr s
\right\}    -\e  . \eqno (3.14) $$ 
Now (3.13) and the definition of $\hat U$ imply the first inequality below, (3.14) the second one while the third one follows from the definition of $\Psi^{n+1}_{c,\l}0$. The equality at the end follows by the definition of $\hat U$. 
$$(\Psi^1_{c,\l}\hat U)(\mu)\ge
E_w\left\{
\int_{0}^{1}
\L_{c,\2\bar\r}(s,\bar X(0,s,\mu),\bar Y)\dr s+\l
\right\}     +$$
$$\liminf_{n\tends+\infty}
E_w\left\{
\int_{1}^{n+1}\L_{c,\2\r_n}(s,X_n(1,s,\r_{\bar Y}(1)),Y_n)\dr s+\l n
\right\}     -\e\ge     $$
$$\liminf_{n\tends+\infty}
E_w\left\{
\int_{0}^{n+1}\L_{c,\2\tilde\r}(s,\tilde X(0,s,\mu),\tilde Y)\dr s  +\l(n+1)
\right\} -2\e \ge
\liminf_{n\tends+\infty}(\Psi^{n+1}_{c,\l}0)(\mu)-2\e  =
\hat U(\mu)-2\e  .  $$
Since $\e>0$ is arbitrary, we get that 
$$(\Psi_{c,\l}^1\hat U)(\mu)\ge\hat U(\mu)\qquad
\forall\mu\in{\cal M}_1(\T^p)  .  \eqno (3.15)$$

On the other hand, let $\mu\in{\cal M}_1(\T^p)$ and let  $Y_n$ minimize in the definition of $(\Psi^{n+1}_{c,\l}0)(\mu)$. Then, by the definition of $\hat U$, we get the equality below.
$$\hat U(\mu)=\liminf_{n\tends+\infty}
E_w\left\{
\int^{n+1}_0\L_{c,\2\r_n}(s,X_n(0,s,\mu),Y_n)\dr s+\l(n+1)
\right\}     $$
where $\r_n$ stays for $\r_{Y_n}$. 
Let $\{ n_h \}$ be a subsequence on which the $\liminf$ is attained. By lemma 3.2, $Y_{n_h}$ is uniformly Lipschitz. Thus, we can apply Ascoli-Arzel\`a\ and, after further refining $\{ n_h \}$, we can suppose that $Y_{n_h}|_{[0,1]\times\T^p}$ converges to a vector field $\bar Y$ in the $C^0$ topology. The formula above yields the first equality below, while the second one follows from the fact that $Y_{n_h}|_{[0,1]\times\T^p}$ converges uniformly to $\bar Y$; the first inequality follows from the definition of 
$\hat U(\hat\r(1)\L^p)$, the second one from the definition of $\Psi^1_{c,0}$.
$$\hat U(\mu)=
\liminf_{h\tends+\infty} E_w\Big\{
\int_0^1\L_{c,\2\r_{n_h}}(s,X_{n_h}(0,s,\mu),Y_{n_h})\dr s+$$
$$\int_1^{n_h+1}\L_{c,\2\r_{n_h}}(s,X_{n_h}(0,s,\r_{n_h}(1)),Y_{n_h})\dr s
+\l(n_h+1)
\Big\}      =   $$
$$E_w\left(
\int_0^1\L_{c,\2\bar\r}(s,\bar X,\bar Y)\dr s  +\l
\right) +$$
$$\liminf_{h\tends+\infty}
E_w\left\{
\int^{n_h+1}_1\L_{c,\2\r_{n_h}}(s,X_{n_h}(1,0,\r_{n_h}(1)),Y_{n_h})\dr s
+\l n_h
\right\}    \ge  $$
$$E_w\left(
\int_0^1\L_{c,\2\bar\r}(s,\bar X,\bar Y)\dr s+\l
\right)   
+\hat U(\bar\r(1)\L^p) \ge
\Psi_{c,0}^1\hat U(\m) \qquad
\forall\mu\in{\cal M}_1(\T^p)   $$
where $\bar\r$ stays for $\r_{\bar Y}$. This proves the inequality opposite to (3.15). Thus, $\hat U=\Psi^1_{c,0}\hat U$; by the last formula, this implies that $\bar Y$ satisfies (5).

It remains to prove that the constant $\l$ is unique. Let 
$\Psi^1_{c,\l_1}\hat U_1=\hat U_1$ and 
$\Psi^2_{c,\l_2}\hat U_2=\hat U_2$. Let us suppose by contradiction that $\l_1<\l_2$. Since $\hat U_i$ is a continuous fixed point, we can suppose that $||\hat U_i||_{\sup}\le M$ for 
$i=1,2$; as a consequence, $\hat U_2\ge\hat U_1-2M$. By proposition 2.10, the first inequality below follows; the second equality follows from the fact that $\hat U_1$ is a fixed point.
$$\Psi^n_{c,\l_2}\hat U_2\ge
\Psi^n_{c,\l_2}\hat U_1-2M\ge
\Psi^n_{c,\l_1}\hat U_1 +n(\l_2-\l_1)-2M=$$
$$\hat U_1 +n(\l_2-\l_1)-2M\ge 
\hat U_2 +n(\l_2-\l_1)-4M  .  $$
For $n$ large enough, the last formula contradicts the fact that 
$\hat U_2$ is a fixed point.

\fin

\vskip 2pc
\centerline{\bf Bibliography}


\noindent [1] L. Ambrosio, Lecture notes on optimal transport problems, mimeographed notes, 2000.

\noindent [2] N. Anantharaman, On the zero-temperature vanishing viscosity limit for certain Markov processes arising from Lagrangian dynamics, J. Eur. Math. Soc. (JEMS), {\bf 6}, 207-276, 2004.

\noindent [3] G. Birkhoff, Lattice theory, vol. {\bf 25}, A. M. S. Colloquium Publ., Providence, R. I., 1967.




\noindent [4] P. Cardialiaguet, Notes on mean field games, from P. L. Lions' lectures at the Coll\`ege de France, 

\noindent mimeographed notes.

\noindent [5] E. Carlen, Lectures on optimal mass transportation and certain of its applications, 2009, mimeographed notes.

\noindent [6] G. Da Prato, Introduction to stochastic analysis and Malliavin calculus, Scuola Normale Superiore, Pisa, 2007.

\noindent [7] R. L. Dobrushin, Vlasov equations, Functional analysis and its applications, {\bf 13}, 45-58, 1979.

\noindent [8] A. Fathi, Weak KAM theorem in Lagrangian dynamics, Fourth preliminary version, mimeographed notes, Lyon, 2003.

\noindent [9] W. H. Fleming, The Cauchy problem for a nonlinear first order PDE, Journal of Differential Equations, {\bf 5}, 515-530, 1969.

\noindent [10] A. Friedman, Partial differential equations of parabolic type, Dover, New York, 1992.

\noindent [11] W. Gangbo, A. Tudorascu, Lagrangian dynamics on an infinite-dimensional torus; a weak KAM theorem, Adv. Math., 
{\bf 224}, 260-292, 2010.

\noindent [12] W. Gangbo, T. Nguyen, A. Tudorascu, Hamilton-Jacobi equations in the Wasserstein space, methods Appl. Anal., {\bf 15}, 155-183, 2008.

\noindent [13] D. Gomes, A stochastic analog of Aubry-Mather theory, Nonlinearity, {\bf 15}, 581-603, 2002.

\noindent [14] T. Hida, Brownian Motion, Springer, Berlin, 1980.




\noindent [15] J. N. Mather, Variational construction of connecting orbits, Ann. Inst. Fourier, {\bf 43}, 1349-1386, 1993.


\noindent [16] M. Viana, Stochastic dynamics of deterministic systems, mimeographed notes, 2000.

\noindent [17] C. Villani, Topics in optimal transportation, Providence, R. I., 2003.

\end

We give a proof of this, following [Tay]. Let $X=C^2(\T^p)$, 
$W=C^1(\T^p)$; we recall that, for $t>0$,
$$||
e^{t\D}
||_{L(W,X)}   \le   \frac{C_{18}}{\sqrt t}  .  \eqno (2.5)$$
For $A>-m$, we let
$${\cal A}_{A}=\{
u\in C((-m,A),X)\st 
\sup_{t\in[-m,A]}||\r-e^{(t+m)\Delta}\d_{x_0}||_{X}\le 1
\}  .  $$
If $\r_1,\r_2\in{\cal A}_{A}$, we set
$$dist(\r_1,\r_2)\colon=\sup_{t\in[-m,A]}||
\r_1(t)-\r_2(t)
||_X    $$
and we note that $dist$ turns ${\cal A}_{A}$ into a complete metric space.
We write $(FP)_{-m,Y,\d_{x_0}}$ as
$$\r(t)=\Psi(\r)(t)\qquad t\in(-1,A)$$
where 
$$\Psi(\r)(t)=e^{(t+m)\D}\d_{x_0}+
\int_{-m}^t e^{(t-s)\D}{\rm div}[\r(s)\cdot Y(s,\cdot)]\dr s  .  $$
We show that, if $A+m$ is small enough, $\Psi$ is a contraction of ${\cal A}_{A}$ into itself. Indeed, for $t\in[-m,A]$, (2.5) implies the second inequality below; the third one follows from the fact that 
$Y\in Vect_m$, and thus $||Y(t,\cdot)||_{C^2(\T^p)}$ is bounded.
$$||
\Phi(\r_1)(t)-\Phi(\r_2)(t)
||_X\le
\int_{-m}^t
||
e^{(t-s)\Delta}{\rm div}[(\r_1(s)-\r_2(s))\cdot Y(s,\cdot)]
||_X  \dr s\le$$
$$\int_{-m}^t\frac{C_{18}}{\sqrt{t-s}}||
{\rm div}[(\r_1(s)-\r_2(s))\cdot Y(s,\cdot)]
||_W\dr s\le
dist(\r_1,\r_2)\cdot\int_{-m}^t\frac{C_{19}}{\sqrt{t-s}}\dr s\le
C_{20}\sqrt{A+m}\cdot dist(\r_1,\r_2)  .  $$
Analogously,
$$||
\Phi(\r)(t)-e^{(t+m)\Delta}\d_{x_0}
||_X\le
\int_{-m}^t
||
e^{(t-s)\Delta}{\rm div}[\r(s)\cdot Y(s,\cdot)]
||_X  \dr s\le$$
$$\int_{-m}^t\frac{C_{21}}{\sqrt{t-s}}\dr s=
C_{22}\sqrt{A+m}\le 1 .  $$
The last two formulas show that $\Phi$ is a contraction of 
${\cal A}_{A}$ into itself, if $A+m$ is small enough. The fixed point $\r_{x_0}$ is the density of a so-called mild solution $\r_{x_0}$ of$(FP)_{-1,Y,\d_{x_0}}$. By the definition of ${\cal A}_A$, we see that $\r_{x_0}$ is $C^2$ in the $x$ variable, and thus that it is a classical solution too. We also note that, by the definition of 
${\cal A}_{A}$, $||\r_{x_0}(A,\cdot)||_{C^2(\T^p)}$ is bounded independently on $x_0$ and $Y\in Vect_m$; thus, standard theorems of PDE theory yield that $\r_{x_0}$ can be prolonged beyond beyond $A$ to 
$t\in(-m,0]$ and that
$$||
\bar\r_{x_0}(t,\cdot)
||_{C([A,0],C^2(\T^p))}\le C_{23}   \eqno (2.6)$$
for a constant $C_{16}$ independent on $x_0$ and on 
$Y\in Vect_m$. On the other side, if $T<A$, we get from the formula above and the definition of ${\cal A}_{A}$ that
$$||\r_{x_0}||_{C^1([T,0]\times\T^p)}+
||\r_{x_0}||_{C((T,0),C^2(\T^p)}\le C_{16}(T)  $$
for a constant $C_{16}(T)$ independent on $x_0\in\T^p$ and on 
$Y\in Vect_m$. Together with (2.6), this implies that the density 
$\r_{x_0}$ satisfies the first formula of (2.3).

Again by the fact that 
$\r_{x_0}\in{\cal A}_{A}$, we get that, as $T\searrow-m$, 
$R_{x_)}(T)=\r_{x_0}\L^p$ converges weak$\ast$ to $\d_{x_0}$; since the Wasserstein metric $d_1$ induces the weak$\ast$ topology, we have that
$$d_1(R_{x_0}(T),\d_{x_0})\le C_{17}(T) \eqno (2.7)$$
with $C_{17}(T)\tends 0$ as $T\tends-m$, independent on $x_0$ and $Y\in Vect_m$.

On the other side,
$$\hat U(\mu)=\liminf_{n\tends+\infty}
E_w\left\{
\int_{-(n+1)}^0\L_{c,\2\r_n}(s,X_n,Y_n)
\right\}   =  $$
$$E_w\left\{
\int_{-(n+1)}^{-n}\L_{c,\2\r_n}(s,X_n,Y_n)
\right\} +
\liminf_{n\tends+\infty}
E_w\left\{
\int_{-(n)}^0\L_{c,\2\r_n}(s,X_n,Y_n)
\right\} .  $$
Taking $Y$ to be the limit of $Y_n$ along a subsequence, we get that $Y$ minimizes and
$$\hat U(\mu)\ge\Lambda_{c,\l}U(\m)  .  $$

We know from proposition 3.1 that $\Lambda^n_{c,\l}U$ is 
$L$-Lipschitz for all $n\in\N$; thus, the thesis follows if we prove that, for some $\l\in\R$,
$$\inf_{k\ge n}\min\mu(\Lambda^k_{c,\l}U)(\mu) \eqno (3.4) $$
is bounded uniformly in $n$. 

Let us consider $\Lambda^n_{c,0}U$; if it remains bounded as 
$n\tends+\infty$, then we have proven the lemma with $\l=0$. Let us suppose now that it is unbounded; to fix ideas, let us suppose that, along a subsequence $n_k$, we have
$$(\Lambda^{n_k}_{c,0}U)\tends-\infty  .  $$
By proposition 3.1, $(\Lambda^{n_k}_{c,0})$ is uniformly Lipschitz on the compact space $\Mt$; thus, the variation
$$\max_\m (\Lambda^{n_k}_{c,0}U)(\m)-
\min_\m (\Lambda^{n_k}_{c,0}U)(\m)$$
is bounded. The last two forumlas imply that we can fix $k$ such that, for some $\a\le\b<0$, we have that
$$\a\le(\Lambda^{n_k}_{c,0}U)\le\b  .  \eqno (3.5) $$
It follows easily from the definition of $\Lambda_{c,0}$ that

\vskip 1pc

\noindent a) if $U_1\le U_2$, then 
$\Lambda_{c,0}U_1\le \Lambda_{c,0}U_2$ and

\noindent b) $\Lambda_{c,0}(U_1+\g)=\Lambda_{c,0}(U_1)+\g$ for all $\g\in\R$. 

\vskip 1pc

Since $\Lambda^n_{c,0}$ is a semigroup, we have that, for 
$s\in\N$ and $0\le l\le n_k-1$, 
$$\Lambda^{sn_k+l}_{c,0}U=(\Lambda^{n_k}_{c,0})^s
\circ\Lambda^l_{c,0}U  $$
so that, setting
$$M=\sup\{
|(\Lambda_{c,0}^lU)(\m)|\st\m\in\Mt,\quad 0\le l\le n_k-1
\}      $$
we get from a), b) and (3.5) that, for $s\in\N$ and $0\le l\le n_k-1$, 
$$s\a-M\le\Lambda_{c,0}^{sn_k+l}U\le s\b+M . \eqno (3.6)$$
Thus, the sequence 
$$\inf_{k\ge n}\left[
\min_\mu\left(
\Lambda^k_{c,0}U
\right)    (\mu)
\right]   $$
grows at most linearly. As a result, 
$$\tilde\l\colon=\liminf_{n\tends+\infty}\frac{1}{n}
\min_\mu\left[
\left(\Lambda^n_{c,0}U
\right)    (\mu)
\right]  $$
is finite. We assert that (3.4) holds for $\l=\tilde\l$.

To show this, we go back to (3.6); since $\Lambda^n_{c,0}U$ is 
$L$-Lipschitz, and $\Mt$ has finite diameter $C_1$, we get that there is $\g_1\in[s\a-M,s\a+M]$ such that
$$\g_1\le\Lambda^{n_k}_{c,0}U\le\g_1+C_1L  .  $$
Iterating as in (3.6), we get that
$$l\g_1\le\Lambda^{ln_k}_{c,0}U\le l(\g_1+C_1L) 
\quad l\ge 1 .  $$
On the other side, again by Lipschitzianity, we have that, for some $\g_2\in[2\g_1,2(\g_1+C_1L)]$,
$$\g_2\le\Lambda^{2n_k}_{c,0}U\le\g_2+C_1L  .  $$
Iterating,
$$l\g_2\le\Lambda^{2ln_k}_{c,0}U\le l(\g_2+C_1L) 
\quad l\ge 1 .  $$
Again by iteration, there is a sequence $\g_r$ such that
$$\g_r\le\Lambda^{2^r n_k}_{c,0}\le (\g_r+C_1L) 
\quad l\ge 1 ,  $$
and
$$\g_l\in[2\g_{l-1},2(\g_{l-1}+C_1L)]  .  $$
Moreover, arguing as in (3.6), we get that
$$a\g_r+s\a-M\le
\Lambda_{c,0}^{a2^r n_k+sn_k+l}U\le
a(\g_r+C_1L)+s\b+M   .   $$
Taking 
$$\tilde\l=\lim_{r\tends+\infty}\frac{\g_r}{2^r}  ,  $$
we get the thesis.

\fin

\vskip 1pc
\noindent{\bf Proof of theorem 4.} We choose $\l$ as in lemma 3.3; thus, we can define
$$\hat U(\m)=\liminf_{n\tends+\infty}
(\Lambda^n_{c,\l}0)(\m)  .  $$
We prove that $\hat U$ is a fixed point. Let $n\in\N$ and let 
$\bar Y$ be a vector field for which
$$E_w\left\{
\int_{-(n+1)}^{-n}\L_{c,\2\bar\r}(s,\bar X,\bar Y)\dr s
\right\}   
+\hat U(\bar\r(-n))\le (\Lambda_{c,\l}\hat U)(\mu)+\e  .   
\eqno (3.8) $$
Let $Y_n$ be such that
$$\min_Y E_w\left\{
\int_{-n}^0\L_{c,\2\r}(s,X,Y)
\right\}  =
E_w\left\{
\int_{-n}^0\L_{c,\2\r_n}(s,X_n,Y_n)
\right\}   .  $$
Let $\tilde Y$ be equal to $\bar Y$ on $[-(n+1),-n]\times\T^p$, to $Y_n$ on $[-n+\d,0]\times\T^p$ and a Lipschitz connection in between. As in proposition 2.9, we can choose the Lipschitz connection in such a way that
$$E_w\left\{
\int_{-(n+1)}^0 
\L_{c,\2\r}(s,\tilde X(-(n+1),s,\mu),\tilde Y(s,\tilde X(-(n+1),s,\mu)))
\dr s
\right\}    \le $$
$$E_w\left\{
\int_{-(n+1)}^{-n} 
\L_{c,\2\r}(s,\bar X(-(n+1),s,\mu),\bar Y(s,\bar X(-(n+1),s,\mu)))
\dr s
\right\}  +  $$
$$E_w\left\{
\int_{-n}^0 
\L_{c,\2\r}(s, X_n(-n,s,\mu),Y_n(s,X_n(-n,s,\mu)))
\dr s
\right\}   +\e  .  $$
Now (3.8) implies the first inequality below, the formula above the second one.
$$(\Lambda_{c,\l}\hat U)(\mu)\ge$$
$$E_w\left\{
\int_{-(n+1)}^{-n}\L_{c,\2\bar\r}(s,\bar X,\bar Y)\dr s
\right\} 
+\liminf_{n\tends+\infty}
E_w\left\{
\int_{-n}^0\L_{c,\2\r_n}(s,X_n,Y_n)
\right\}     -\e\ge$$
$$\liminf_{n\tends+\infty}
E_w\left\{
\int_{-(n+1)}^{0}\L_{c,\2\tilde\r}(s,\tilde X,\tilde Y)\dr s
\right\} -2\e\ge \hat U(\mu)-2\e  .  $$
On the other side,
$$\hat U(\mu)=\liminf_{n\tends+\infty}
E_w\left\{
\int_{-(n+1)}^0\L_{c,\2\r_n}(s,X_n,Y_n)
\right\}   =  $$
$$E_w\left\{
\int_{-(n+1)}^{-n}\L_{c,\2\r_n}(s,X_n,Y_n)
\right\} +
\liminf_{n\tends+\infty}
E_w\left\{
\int_{-(n)}^0\L_{c,\2\r_n}(s,X_n,Y_n)
\right\} .  $$
Taking $Y$ to be the limit of $Y_n$ along a subsequence, we get that $Y$ minimizes and
$$\hat U(\mu)\ge\Lambda_{c,\l}U(\m)  .  $$